\documentclass[11pt]{article}

\usepackage[T1]{fontenc} 
\usepackage[latin1]{inputenc}

\usepackage{hyperref}
  
\usepackage[english]{babel}
\usepackage{xy}
\usepackage[T1]{fontenc} 
\usepackage[latin1]{inputenc}
\usepackage{amsfonts,amssymb,amsmath}
\usepackage{hyperref} 
\usepackage{graphics} 
\usepackage{multido, xcolor}
 
\newtheorem{theo}{\bf Theorem}[section]
\renewcommand{\theequation}{\thesection.\arabic{equation}}
\renewcommand{\thesubsection}{\thesection.\arabic{subsection}}
\renewcommand{\theequation}{
\arabic{equation}}

\newcommand{\Section}{\setcounter{equation}{0} \section}

\newtheorem{coro}[theo]{\sc Corollary}
\newtheorem{defi}[theo]{\sc Definition}

\newtheorem{prop}[theo]{\sc Proposition}

\newtheorem{main}[theo]{\sc Main Theorem}

\newtheorem{exam}[theo]{\sc Example}

\def\qed{ \ \hfil$\square$}

\def\Ker{{\hskip0.3mm\rm  Ker\hskip0.5mm}}
\def\Im{{\hskip0.3mm\rm  Im\hskip0.5mm}}

\makeindex

\def \F {{\mathbb F}}
\def\Si{\Sigma}

\def\ph{\varphi}
\def\W{A}
\def\End{{\rm End}}
\def\is{\buildrel \sim   \over \rightarrow}

\newcommand{\di}[1]{\displaystyle{#1}}

\newcommand{\tr}{{\rm tr}}

\newcommand{\dd}{{\rm d}}

\newcommand{\Gal}{{\rm Gal}}

\newcommand{\GL}{{\rm GL}}
\newcommand{\GSp}{{\rm GSp}}
\newcommand{\Hol}{{\mathcal Hol}}
\newcommand{\Hom}{{\rm Hom}}

\newcommand{\M}{{{\rm M}}}

\newcommand{\ord}{{\rm ord}}
\newcommand{\SL}{{\rm SL}}
\newcommand{\SO}{{\rm SO}}

\newcommand{\Tr}{{\rm Tr}}

\newcommand{\A}{{\mathbb A}}
\newcommand{\C}{{\mathbb C}}
\newcommand{\N}{{\mathbb N}}
\newcommand{\Q}{{\mathbb Q}}
\newcommand{\R}{{\mathbb R}}
\newcommand{\Z}{{\mathbb Z}}
\newcommand{\Hb}{{\mathbb H}}

\newcommand{\Qb}{\overline{\mathbb Q}}

\newcommand{\chib}{\bar{\chi}}

\newcommand{\Ac}{{\mathcal A}}
\newcommand{\Br}{{\mathcal B}}

\newcommand{\Cc}{{\mathcal C}}
\newcommand{\Ci}{{\mathcal C}^{\infty}}

\newcommand{\Mc}{{\mathcal M}}

\newcommand{\Or}{{\mathcal O}}
\newcommand{\Oc}{{\mathcal O}}

\newcommand{\Sc}{{\mathcal S}}
\newcommand{\Pc}{{\mathcal P}}

\newcommand{\X}{{\mathcal X}}

\newcommand{\al}{\alpha}
\newcommand{\De}{\Delta}
\newcommand{\de}{\delta}
\newcommand{\e}{\varepsilon}
\newcommand{\ep}{\varepsilon}
\newcommand{\Ga}{\Gamma}
\newcommand{\ga}{\gamma}

\newcommand{\ka}{\kappa}
\newcommand{\La}{\Lambda}
\newcommand{\la}{\lambda}

\font\teneusm=eusm10 \font\seveneusm=eusm7 
\font\fiveeusm=eusm5 
\newfam\eusmfam 
\textfont\eusmfam=\teneusm 
\scriptfont\eusmfam=\seveneusm 
\scriptscriptfont\eusmfam=\fiveeusm

\def\nd{\not\hskip3.0pt\mid}
\def\nds{\not\hskip1.5pt\mid}
\def\mat #1,#2,#3,#4,{\left({#1\atop #3}{#2\atop #4}\right)}
\def\bra#1,{{\left\lbrace {#1}\right\rbrace}}
\def\ph{\varphi}
\def\Ph{\Phi}
\def\si{\sigma}
\def\z{\zeta}
\def \sign {{\rm sign}}
\def \sgn {{\rm sgn}}
\def \inv {{\rm inv}}
\def \diag {{\rm diag}}
\def \Inf {{\rm Inf}}
\def \Proj {{\rm Proj}}
\def \Mor {{\rm Mor}}
\def \Div {{\rm Div}}
\def \div {{\rm div}}
\def \Pic {{\rm Pic}}
\def \Res {{\rm Res}}
\def \vol {{\rm vol}}
\def \Ver {{\rm Ver}}
\def \Spin {{\rm Spin}}
\def \Spec {{\rm Spec}}
\def \Step {{\rm Step}}
\def \Spl {{\rm Spl}}
\def \M{{\rm M}}
\def \Gal{{\rm Gal}}
\def \suc{{\rm suc}}
\def \Gab{G^{\rm ab}}
\def \rk{{\rm rk}}
\def \Supp{{\rm Supp}}
\def \Id{{\rm Id}}
\def \Sp{{\rm Sp}}
\def\tm{\tilde m}
\def\tGa{\tilde{\Ga}}
\def\is{\buildrel \sim   \over \rightarrow}
\def\bs{\backslash}
\def\rr{\rangle}
\def\l1{\langle}

\newcommand{\B}{\left(\begin{array}{cc}}
\newcommand{\E}{\end{array}\right)}
\newcommand{\p}{\tilde}

\def \fns{{${}^{*)}$}}

\newcommand{\comm}[1]
{\fns\marginpar{$\boxed
{\hskip-6pt
{\small {\sf 
\begin{tabular} {l}
 #1
\end{tabular}
}
}
}
$
}
}
\def \?  {\comm{check ?}}

\usepackage{euscript}
\let\scr=\EuScript

\def\Hol{{\scr Hol}}
\let\mathcal=\scr           \def\Or{{\scr O}}
\def\ang#1,{{\left\langle {#1}\right\rangle}} 

\def\G{{\dbG}}

\def \rem  {\medskip\noindent {\sc Remarque}}

\def\qed{ \ \hfil$\square$}

\def\Ker{{\hskip0.3mm\rm  Ker\hskip0.5mm}}
\def\Im{{\hskip0.3mm\rm  Im\hskip0.5mm}}

\def \F {{\mathbb F}}
\def\Si{\Sigma}

\def\ph{\varphi}
\def\W{A}
\def\End{{\rm End}}
\def\is{\buildrel \sim   \over \rightarrow}

\font\teneusm=eusm10 \font\seveneusm=eusm7 
\font\fiveeusm=eusm5 
\newfam\eusmfam 
\textfont\eusmfam=\teneusm 
\scriptfont\eusmfam=\seveneusm 
\scriptscriptfont\eusmfam=\fiveeusm

\font\tengothic=eufm10
\font\sevengothic=eufm7
\font\fivegothic=eufm5
\newfam\Gothic
\textfont\Gothic=\tengothic
\scriptfont\Gothic=\sevengothic
\scriptscriptfont\Gothic=\fivegothic
\def\go{\fam\Gothic\tengothic}

\def\nd{\not\hskip3.0pt\mid}
\def\nds{\not\hskip1.5pt\mid}
\def\mat #1,#2,#3,#4,{\left({#1\atop #3}{#2\atop #4}\right)}
\def\bra#1,{{\left\lbrace {#1}\right\rbrace}}
\def\ph{\varphi}
\def\Ph{\Phi}
\def\si{\sigma}
\def\z{\zeta}
\def\w{\omega}

\def\W{\Omega}
\def \sign {{\rm sign}}
\def \sgn {{\rm sgn}}
\def \inv {{\rm inv}}
\def \diag {{\rm diag}}
\def \Inf {{\rm Inf}}
\def \Proj {{\rm Proj}}
\def \Mor {{\rm Mor}}
\def \Div {{\rm Div}}
\def \div {{\rm div}}
\def \Pic {{\rm Pic}}
\def \Res {{\rm Res}}
\def \Re {{\rm Re}}
\def \vol {{\rm vol}}
\def \Ver {{\rm Ver}}
\def \Spin {{\rm Spin}}
\def \Card {{\rm Card}\,}
\def \Step {{\rm Step}}
\def \Spl {{\rm Spl}}
\def \M{{\rm M}}
\def \Gal{{\rm Gal}}
\def \suc{{\rm suc}}
\def \Gab{G^{\rm ab}}
\def \rk{{\rm rk}}
\def \Supp{{\rm Supp}}
\def \Id{{\rm Id}}
\def \Sp{{\rm Sp}}
\def\m{{\go m}}
\def\tGa{\tilde{\Ga}}
\def\is{\buildrel \sim   \over \rightarrow}
\def\bs{\backslash}
\def\rr{\rangle}
\def\l1{\langle}
\def\Eq{\Longleftrightarrow}
\def\Numero{${\rm N\sp\circ}$}

\let\scr=\EuScript

\def\Hol{{\scr Hol}}
\let\mathcal=\scr           \def\Or{{\scr O}}

\let\bold=\boldsymbol

\def\G{{\dbG}}

\def\vin{{ {\tiny \mid }  
\kern-7.29pt 
\bigcup }}

\def \bs{{\backslash}}

\def\ang#1,{{\left\langle {#1}\right\rangle}}

\def \Aut {{\rm Aut}}
\def\m{{\frak m}}

\normalsize

\def\ab{\A}

\def\qb{\Q}
\def\rb{\R}
\def\zb{\Z}

\def\la{\displaystyle\mathop{\longrightarrow}}

\def\dc{\Dr}
\def\hc{\Hr}

\def\mc{\Mr}

           \def\ab{{\Bbb A}}            
                      
           \def\cb{{\Bbb C}}            
\def\dc{{\cal D}}           

\def\hc{{\cal H}}

\def\mc{{\cal M}}

           \def\qb{{\Bbb Q}}            
           \def\rb{{\Bbb R}}

           \def\zb{{\Bbb Z}}

\def\diag {\mathop{\rm diag}\nolimits}
\def\SL {\mathop{\rm SL}\nolimits}
\def\Sp {\mathop{\rm Sp}\nolimits}

\font\tenrus=wncyr10
\font\sevenrus=wncyr7
\newfam\Rus
\textfont\Rus=\tenrus
\scriptfont\Rus=\sevenrus
\def\ru{\fam\Rus\tenrus}
\def\scr{\fam\Rus\sevenrus} 

\font\tenbrus=wncyb10
\font\eightbrus=wncyb8
\newfam\Brus
\textfont\Brus=\tenbrus
\scriptfont\Brus=\eightbrus

\font\tenirus=wncyi10
\font\eightirus=wncyi8
\newfam\Irus
\textfont\Irus=\tenirus
\scriptfont\Irus=\eightirus

\font\tenrus=wncyr10
\font\sevenrus=wncyr7
\newfam\Rus
\textfont\Rus=\tenrus
\scriptfont\Rus=\sevenrus
\def\ru{\fam\Rus\tenrus}
\def\scr{\fam\Rus\sevenrus} 
\input cyracc.def
\def\i1{\accent'044i}
\def\I1{\accent'044I}
\def\e1{\accent'040e}
\def\l1{l{}p1}

\def\qed{\quad\hbox{\hskip 1pt\vrule width 4pt height 6pt
          depth 1.5pt\hskip 1pt}}

\title
{
Analytic constructions of
$p$-adic $L$-functions and Eisenstein series
}

\author{Alexei PANCHISHKIN\\ \small
Institut Fourier, 
  Université Grenoble-1\\
B.P.74, 38402 St.--Martin d'H\`eres, FRANCE\\
\it  Conference 
Automorphic Forms and Related Geometry, \\
\it Assessing the Legacy of I.I. Piatetski-Shapiro}

\date
{23 - 27 April, 2012 (Yale University in New Haven, CT)}

\def\Numero{${\rm N\sp\circ}$}
\newcounter{ncours}{\setcounter{ncours} {1}}

\begin{document}

\maketitle
\subsection*{Work of Ilya Piatetski-Shapiro and new ways \\ of constructing 
complex and $p$-adic 
$L$-functions}

In June 2011 Roger Howe invited me to 
this conference
devoted to assessing the work of Ilya Piatetski-Shapiro, especially
in the areas of autmorphic forms and geometry. 

\

Many thanks to the organizers for this invitation and this occasion 
both to review the accomplishments of Ilya Piatetski-Shapiro and his
colleagues, 
and to point to productive directions to take research
from here.

\

I new Ilya since 1973-74 during our joint participation
 in the seminar of Manin and Kirillov on p-adic L-functions,
and attending his informal lectures on GL(3) in Moscow University in
April-May 1975.

\

After many years we met again  in Jerusalem in February 1998 during the conference "p-Adic Aspects of the Theory of Automorphic Representations". 

\

Ilya liked my construction of $p$-adic standard $L$-functions of Siegel modular forms \cite{Pa91}, and suggested to extend it to spinor $L$-functions, using the restriction of an Eisenstein series
 to the Bessel subgroup in the generalized Whittaker models (see Olga Taussky Todd memorial volume \cite{PS3}).
So we started a joint work "On $p$-adic $L$-functions for $GSp(4)$".

\

\subsection*{
In 1998,  a conference for Ilya Piatetski-Shapiro} was organized in the Fourier Institute
(Grenoble, France),
with participation of A.Andrianov, G.Henniart, H.Hida, J.-P.Labesse,
J.-L.Waldspurger and others.

\

In IAS, we had the most intensive period of our joint work in the Fall 
1999-2000.

My last meeting Ilya was  on January 31,
2009 
at
the  the Weizmann Institute home of 
Volodya  Berkovich and his wife Lena, who is
Grisha Freiman's daughter.
Vera brought Ilya,
Edith, and Edith's mother Ida to Berkovich's.
This
very pleasant gathering also included Grisha Freiman, his wife Nina, 
Steve and Mary Gelbart, Antoine Ducros, and my wife Marina, 
see also \cite{CGS}, p. 1268.

\

That was the last party in Ilya's life.

\

His favorite automorphic forms were the {\it Eisenstein series}, 
and the main subject of this talk will be my new  construction 
of {\it meromorphic $p$-adic families of Siegel-Eisenstein series
in relation to the geometry of homogeneuos spaces}, both complex and $p$-adic, for any prime $p$.

I am glad that this construction fits into the particular subject
"Automorphic Forms and Related Geometry" of our conference.

\

\subsection*{ $p$-adic Siegel-Eisenstein series and related geometry}
Let us consider the {\it symplectic group $\Gamma=\Sp_{m}(\Z)$} (of $(2m\times 2m)$-matrices), and prove that
the {\it Fourier coefficients $a_h(k)$ of the original Siegel-Eisenstein
series $E_k^m$} admit an explicit $p$-adic meromorphic interpolation on $k$
where $h$ runs through all positive definite half integral matricies
for $\det(2h)$ not divisible by $p$, where
\begin{align*}
E_k^{m}(z)&=\sum_{(c,d)/\sim}\det(cz+d)^{-k}=\sum_{\gamma\in P\backslash \Gamma}\det(cz+d)^{-k}
\\ &=\sum_{h\in B_m}a_h\exp(\tr (hz))
\end{align*}
on the Siegel upper half plane  
${\mathfrak H}_m=\{z={}^t z\in M_m(\C) | \Im (z)>0\}$ of degree $m$,  $(c,d)$ runs over equivalence classes of all coprime 
symmetric couples, $\gamma= \begin{pmatrix} a&b\\c&d\end{pmatrix}$ runs over equivalence classes
of $\Gamma$ 
modulo the Siegel parabolic $P=\begin{pmatrix} *&*\\0&* \end{pmatrix}$.

\

\subsection*{ $p$-adic Siegel-Eisenstein series and related geometry}

{\it The homogeneuos  space $X=\{(c,d)/\sim\}=P\backslash \Sp_{m} $ and its $p$-adic points admit Siegel's  coordinates} 
$$
\nu=\det (c)\mbox{ and }
{\mathfrak R}
=c^{-1}d
$$ 
defined on the main subset  given by $\det (c)\in \mathrm{GL}_1$,
which is  used in the construction.

\

I try also to present  various {\it applications:  to  $p$-adic $L$-functions, 
to Siegel's Mass Formula, to $p$-adic analytic families of automorphic representations}.

\


 {\it Eisenstein series} are basic automorphic forms,  
and there exist several ways to construct them via group theory, lattice theory, Galois representations, spectral theory...


\

For me, 
the Eisenstein series is the main tool of analytic constructions of
 complex and $p$-adic $L$-functions, in particular via the {\it doubling method}, 
see \cite{PSR}, \cite{GRPS}, \cite{Boe85}, \cite{Shi95}, \cite{Boe-Sch},\dots,   
greatly thanks to  Ilya Piatetski-Shapiro and his collaborators.



\

\subsection*{General strategy
 }
 For any Dirichlet character $\chi\bmod p^v$ 
 consider Shimura's "involuted" Siegel-Eisenstein series
 assuming their absolute convergence (i.e. $k> m+1$):
 $$
E^*_k(\chi,z)=\sum_{(c,d)/\sim}\chi(\det(c))\det(cz+d)^{-k}
=\sum_{0<h\in B_m}a_h(k, \chi)q^h.
$$
The two sides of the equality produce {\it dual approaches: geometric and algebraic}.
The Fourier coefficients can be computed by Siegel's method (see \cite{St81}, 
\cite{Shi95}, \dots )
via the singular series
\begin{align}\label{ahk}
&a_h({E}^*_k(\chi,z))\\ &\nonumber = 
\frac {(-2\pi i)^{mk}} {2^{\frac{m(m-1)} {2}}\Gamma_m(k)} 
\sum_{{\mathfrak R}\bmod 1}\chi(\nu({\mathfrak R}))\nu({\mathfrak R})^{-k}
\det h^{k-\frac{m+1} 2}e_m(h{\mathfrak R})
\end{align}
The {\it orthogonality relations $\bmod p^v$ produce two families of distributions}
(notice that terms in the RHS are invariant under sign changes, and (\ref{ah0}) is algebraic after multiplying by the factor in (\ref{ahk})):
\small
\begin{align}&\label{st}
\hskip-1.6cm{1\over \varphi(p^v)} \sum_{\chi\bmod p^v} \bar\chi (b)
 \sum_{(c,d)/\sim}{\chi(\det(c))\over \det(cz+d)^{k}}=
 \sum_{(c,d)/\sim\atop 
\det(c)\equiv b \bmod p^v}{\sgn(\det(c))^k\over \det(cz+d)^{k}} 
\end{align}
\hskip-1.3cm\begin{align} 
\label{ah0} & 
{
 1\over \varphi(p^v)} \sum_{\chi\bmod p^v} \bar\chi (b)
\sum_{{\mathfrak R}\bmod 1}{\chi(\nu({\mathfrak R}))e_m(h{\mathfrak R})\over \nu({\mathfrak R})^{k}}
=
\sum_{{\mathfrak R}\bmod 1\atop 
\nu({\mathfrak R})\equiv b \bmod p^v}{e_m(h{\mathfrak R})\sgn\nu({\mathfrak R})^k\over \nu({\mathfrak R})^{k}}
\end{align}
\subsection*{The use of Iwasawa theory and  pseudomeasures}

We express the integrals of Dirichlet characters $\theta \bmod p$
along the distributions (\ref{ah0})
through the reciprocal of a product of $L$-functions, and  elementary integral factors.
The result turns out to be an {\it Iwasawa function of the variable $t=(1+p)^k-1$ divided by a
distinguished polynomial} provided that $\det h$ is not divisible by $p$.

Thus the second family (\ref{ah0})
comes from a  unique pseudomeasure $\mu_h^*$  which becomes a measure after multiplication by an explicit polynomial factor (in the sense of the convolution product).

\

Then we deduce that
(\ref{st}) determines a unique pseudomeasure with coefficients in 
${\Q}[\![q^{B_m}]\!]$ whose moments are given by those of the coefficients (\ref{ah0})
(after removing from the Fourier expansion 
$f(z)=\sum_{h\ge 0}a_he_m(hz)$
all $h$ with $\det h$ divisible by $p$):
$$
\sum_{h>0, p\nds \det h}a_hq^h
=
p^{-m(m+1)/2}
\sum_{h_0\bmod p\atop p\nds \det h_0}\sum_{x\in S\bmod p}e_m(-h_0x/p)f(z+(x/p)).
$$
\

In this way a $p$-adic family of Siegel-Eisenstein series is geometrically produced.

\


\tableofcontents

\section{Complex and $p$-adic $L$-functions}
\subsection*{Generalities about $p$-adic $L$-functions}
There exist two kinds of $L$-functions
\begin{itemize}

\item Complex-analytic $L$-functions (Euler products) 
\item $p$-adic $L$-functions (Mellin transforms $L_\mu $ of $p$-adic measures) 

\end{itemize}
Both are used in order to obtain a number ($L$-value) from an automorphic form.
Usually such a number is algebraic (after normalization) via the embeddings 
$$
\overline \Q\hookrightarrow\C, \ \ \overline \Q\hookrightarrow\C_p=\widehat{\overline\Q}_p.
$$
{\it How to define and to compute $p$-adic $L$-functions?} 
We use Mellin transform of a $\Z_p$-valued distribution $\mu$ on a profinite group
$$
Y=
\lim_{\buildrel\leftarrow\over i} Y_i, \ \mu\in {\it Distr}(Y, \Z_p)=\Z_p[[Y]] =
\lim_{\buildrel\leftarrow\over i}\Z_p [Y_i]=:\Lambda_Y 
$$
(the {\it Iwasawa algebra} of  $Y$).
$$
L_\mu (x)=\int_Yx(y)d \mu, \ \ x\in X_Y=\Hom_{cont}(Y, \C_p^*) 
$$
(the {\it Mellin transform } of  $\mu$ on $Y$).

\

\subsection*{Examples of p-adic measures and $L$-functions} 
\begin{itemize}
\item $Y=\Z_p$, $X_Y=\{\chi_t: y\mapsto (1+t)^y\}$. 
The Mellin transform  
$$L_\mu(\chi_t)=\int_{\Z_p} (1+t)^y d\mu (y)$$ of any  measure $\mu$
on $\Z_p$ is  given by {\it the Amice transform},
which is the following power series
$$ 
A_\mu(t)=\sum_{n\ge 0}t^n\int_{\Z_p} {y\choose n} d\mu (y)=\int_{\Z_p} (1+t)^y d\mu (y), 
$$
e.g. $A_{\delta_m}=(1+t)^m$. Thus,
${\it Distr}(\Z_p, \Z_p)\cong\Z_p[\![T]\!]$.
\item
$
Y=\Z_p^*=\Delta \times \Gamma= \{y=\delta (1+p)^z, \delta^{p-1}=1, z\in\Z_p\}
$ \\
$X_{\Z_p^*}
=\{\theta \chi_{(t)}\ |\ \theta \bmod p, \chi_{(t)}(\chi_{(t)})=(1+t)^z$,  
where \\
$\Delta$ is the  subgroup of  roots of unity,
$\Gamma=1+p\Z_p$.\\
The {\it $p$-adic Mellin transform}
$L_\mu(\theta \chi_{(t)})=\int_{\Z_p^*}\theta(\delta)(1+t)^z \mu(y)$
 of a measure $\mu$ on  $\Z_p^*$ 
 is given by the collection of Iwasawa series 
$\displaystyle
G_{\theta, \mu} (t)=  
\sum_{n\ge 0} a_{n,\theta} t^n$, where $\displaystyle (1+t)^z=\sum_{n\ge 0}{z\choose n}t^n$,
\\
$\displaystyle
 a_{n,\theta}=\sum_{\delta \bmod p,\  n\ge 0}\theta (\delta)t^n\cdot
\int_{\Z_p} {z\choose n} \mu (\delta (1+p)^z) .
$



\item
A general idea is to construct
$p$-adic   $L$-functions  {\it directly from  Fourier coefficients} of modular forms (or from the Whittaker functions of automorphic forms).

\end{itemize}

\


\section{$p$-adic meromorphic continuation of the Siegel-Eisenstein series}

\subsection*{Mazur's $p$-adic integral 
}
For any choice of a natural number $c\ge 1$ not divisible by $p$, there exists a $p$-adic measure 
$\mu_c$ on $\Z_p^*$, such that the special values
$$
\zeta(1-k)(1-p^{k-1})=
{\int_{\Z_p^*}y^{k-1}d\mu_c \over 1-c^{k}}
\in \Q \ \ (k \ge 2\mbox{ even }),
$$
produce the {\it Kubota-Leopoldt $p$-adic zeta-function} 
$
\zeta_p: X_p \to \cb_p
$ 
(where $X_p=X_{\Z^*_p}=Hom_{cont}(\Z_p^*, \cb_p^*))$ 
as the  {\it $p$-adic Mellin transform}
$$
\zeta_p(x)=
{\int_{\Z_p^*}x(y)d\mu_c (y)\over 1-cx(c)}= {L_{\mu_c}(x)\over 1-cx(c)}, 
$$
with a single simple pole at 
$x=x_p^{-1}\in X$, 
where $\cb_p=\hat{\overline{\qb}}_p$ the Tate field, the completion of an algebraic closure of the  $p$-adic field
 ${\qb}_p$,  $x\in X_p$ (a $\cb_p$-analytic Lie group), 
$x_p(y)=y\in X_p$, and $x(y)=\chi(y)y^{k-1}$ as above.

{\large \it Explicitly}: Mazur's measure is given by 
$$
\mu_c(a+p^v\Z_p)=
\frac {1}c\left[ \frac {ca}{p^v}\right] + 
\frac {1-c}{2c}=\frac{1}c B_1(\{\frac {ca}{p^v}\})-B_1(\frac {a}{p^v}),~B_1(x)=x-\frac 1 2,
$$
see \cite{LangMF}, Ch.XIII.

\

\subsection*{Meromorphic $p$-adic continuation  of ${1\over \zeta(1-k) (1-p^{k-1})}$} 

 For any odd prime $p$ take the Iwasawa series $G_{\theta,c}(t)$ of Mazur's measure $\mu_c$ where $\theta$ is a character $\bmod p$, 
$\displaystyle G_{\theta,c}(t):=\int_{\Z_p^*}\theta(y)\chi_{(t)}(\langle y\rangle)\mu_c
=\sum_{n=0}^\infty a_nt^n\in \Z_p[\![t]\!], 
$~and 
\\ $\displaystyle
\chi_{(t)}:(1+p)^z\mapsto (1+t)^z,  \  \langle y\rangle=\frac y {\omega(y)}
, \omega$  the Teichmüller character.
Mazur's integral of the character  $y^{k-1}=\omega^{k-1}\cdot \chi_{(t)}$ shows that $\theta=\omega^{k-1}, (1+t)=(1+p)^{k-1}$ 
\begin{align}\label{Gthc}
\zeta(1-k)(1-p^{k-1})
={G_{\theta,c}((1+p)^{k-1}-1)  \over 1-c^k}.
\end{align}
By the {\it Weierstrass preparation theorem}  we have a decomposition 
$$
G_{\theta,c}(t)=U_{\theta,c}(t) P_{\theta,c}(t)
$$
with a {\it distinguished polynomial} $P_{\theta,c}(t)$ and {\it invertible power series} $U_{\theta,c}(t)$. The inversion of (\ref{Gthc}) for any even $k\ge 2$ gives :

\begin{align*}
{1\over \zeta(1-k)(1-p^{k-1})}
=G_{\theta,c}((1+p)^{k-1}-1)^{-1}  (1-c^k).
\end{align*}

\subsection*{The answer: for any prime $p>2$ and  even $k\ge 2$}  
is the following Iwasawa function on $t=t_k=(1+p)^{k}-1$ divided by a distingushed polynomial:
\begin{align}\label{zki}
\frac 1 {\zeta(1-k)(1-p^{k-1})}&=\frac {U^*_{\theta,c}((1+p)^{k-1}-1) (1-c^k)}
{P_{\theta,c}((1+p)^{k-1}-1)}\\ & \nonumber
=\frac {U^*_{\theta,c}((1+t_k)(1+p)^{-1}-1) (1-c^k)}
{P_{\theta,c}((1+t_k)(1+p)^{-1}-1)}
\end{align}
which is meromorphic in the unit disc of the variable $t=(1+p)^k-1$
with a finite number of poles (expressed via roots of $P_{\theta,c}$) for $\theta=\omega^{k-1}$,
and
$$
U^*_{\theta,c}((1+p)^{k-1}-1):=1/U_{\theta,c}((1+p)^{k-1}-1).
$$

The above formula immediately extends to all Dirichlet $L$-functions
of characters $\chi \bmod p^v$ as the following Iwasawa function  divided by a polynomial:
$$
\frac 1 {L(1-k, \chi)(1-\chi(p)p^{k-1})}
= \frac {U^*_{\theta,c}(\chi(1+p)(1+p)^{k-1}-1) (1-\chi(c)c^k)}
{P_{\theta,c}(\chi(1+p) (1+p)^{k-1}-1)}
$$
where 
$
U^*_{\theta,c}(\chi(1+p)(1+p)^{k-1}-1):={\displaystyle 1\over \displaystyle U_{\theta,c}(\chi(1+p)(1+p)^{k-1}-1)}
$ 

\

\subsection*{Illustration: numerical values of $\zeta(1-2k)^{-1}(1-p^{2k-1})^{-1}$ for $p=37$}  
gp >{\tt zetap1(p,n)= -2*n/(bernfrac(2*n)*(1-{p}\^{}{(2*n-1)}+O({p}\^{}5)));}\\
gp > p=37;\\
gp > for(k=1,({p}-1)/2, print(2*k,  zetap1(p,k)))\\
\begin{tabular}{l|l}\hline
$2k$& $\zeta(1-2k)^{-1}(1-p^{2k-1})^{-1}$ \\
 \hline
2       &$25 + 24*37 + 24*37^2 + 24*37^3 + 24*37^4 + O(37^5)$\\ \hline
4       &$9 + 3*37 + 9*37^3 + 3*37^4 + O(37^5)$\\ \hline
6       &$7 + 30*37 + 36*37^2 + 36*37^3 + 36*37^4 + O(37^5)$\\ \hline
8       &$18 + 6*37 + O(37^5)$\\ \hline
10       &$16 + 33*37 + 36*37^2 + 36*37^3 + 36*37^4 + O(37^5)$\\ \hline
12       &$8 + 25*37 + 28*37^2 + 23*37^3 + O(37^5)$\\ \hline
14       &$25 + 36*37 + 36*37^2 + 36*37^3 + 36*37^4 + O(37^5)$\\ \hline
16       &$6 + 16*37 + 31*37^2 + 29*37^3 + 20*37^4 + O(37^5)$\\ \hline
18       &$3 + 4*37 + 10*37^2 + 32*37^3 + 25*37^4 + O(37^5)$\\ \hline
20      &$11 + 13*37 + 19*37^2 + 36*37^3 + 12*37^4 + O(37^5)$\\ \hline
22      &$1 + 26*37 + 15*37^2 + 35*37^3 + 9*37^4 + O(37^5)$\\ \hline
24      &$16 + 28*37 + 24*37^2 + 27*37^3 + 31*37^4 + O(37^5)$\\ \hline
26      &$4 + 17*37 + 25*37^2 + 25*37^3 + 19*37^4 + O(37^5)$\\ \hline
28      &$22 + 36*37 + 8*37^2 + 4*37^3 + 33*37^4 + O(37^5)$\\ \hline
30      &$22 + 5*37 + 35*37^2 + 9*37^3 + 5*37^4 + O(37^5)$\\ \hline
{\bf 32}      &$36*37^{-1} + 28 + 3*37 + 19*37^2 + 18*37^3 + O(37^4)$\\ \hline
34      &$20 + 37 + 30*37^2 + 15*37^3 + 22*37^4 + O(37^5)$\\ \hline
\fbox{36}      &$36*37 + 29*37^2 + 35*37^3 + 5*37^4 + 37^5 + O(37^6)$\\ 
\hline
\end{tabular}
\\

\subsection*{Fourier expansion of the Siegel-Eisenstein series}
has the form
$$
E_k^{m}(z)=\sum_{\gamma\in P\backslash \Gamma}\det(cz+d)^{-k}
=\sum_{h\in B_m}a_hq^h, 
$$
where $a_h=a_h(k)=a_h(E_k^{m})$, $q^h=e^{2\pi i\tr (hz)}$, $h$ runs over semi-definite half integral $m\times m$ matrices.

\

The {\it rationality of the coefficients $a_h$} 
was established in Siegel's pioneer work \cite{Si35}  in connection with a study of local densities for quadratic forms. 
Siegel expressed
$a_h(k)$ as a product of local factors over all primes and $\infty$. 

In a difficult later work \cite{Si64b} Siegel proved the boundedness of their denominators, 
and S.Boecherer \cite{Boe84} gave a simplified proof of a more precise result in 1984.
M.Harris extended the rationality to  
wide classes of Eisenstein series on Shimura varieties
\cite{Ha81}, \cite{Ha84}. Their relation to the Iwasawa Main Conjecture
and $p$-adic $L$-functions on the unitary groups
was established in   \cite{HLiSk}.

\

\subsection*{Explicit $p$-adic continuation of $a_h(k)$}
as  {\it Iwasawa functions on $t=(1+p)^{k}-1$ divided by distinguished polynomials}.
Let $a_h^{(p)}(k)$ denote the {\it $p$-regular part of the coefficient $a_h(k)$} 
(i.e. with the Euler $p$-factor removed from the product). 
Namely, for any even $k$, $a_h^{(p)}(k)=a_h(E_k^{m})$ times 
$$
\begin{cases} 1/((1-p^{k-1})(1+\psi_h (p)p^{k-\frac m 2-1})\prod_{i=1}^{(m/2)-1}(1-p^{2k-2i-1})) & \\
= (1-\psi_h (p)p^{k-\frac m 2-1})/((1-p^{k-1})\prod_{i=1}^{m/2}(1-p^{2k-2i-1})),& m \mbox{ even}
\\ 
1/((1-p^{k-1})\prod_{i=1}^{(m-1)/2}(1-p^{2k-2i-1})),& m \mbox{ odd},  
\end{cases}
$$ 
where the {\it $p$-correcting factor is a $p$-adic unit}, 
and $\displaystyle \psi_h (n):=\left({\det(2h)(-1)^{m/2} \over n}\right)$.

\begin{main}[
A.P., 2012]\label{pMSE}
Let $h$ be any positive definite half integral matrix
with $\det(2h)$ not divisible by $p$.
Then there exist explicitly given distinguished polynomials $P^E_{\theta,h}(T)\in \Z_p[T]$
and Iwasawa series $S^E_{\theta,h}(T) \in \Z_p[\![T]\!]$
such that the $p$-regular part $a_h^{(p)}(k)$ of the Fourier coefficient $a_h(k)$ 
admit the following $p$-adic meromorphic interpolation on all even $k$ with $\theta=\omega^k$ fixed 
$$
a_h^{(p)}(k)=
{S^E_{\theta,h}((1+p)^{k-1}-1)\over P^E_{\theta,h}((1+p)^{k-1}-1)}
$$
with a finite number of poles expressed via the roots of  $P^E_{\theta,h}(T)$
where the denominator depends only on 
$\det(2h) \bmod 4p$ and $k   \bmod p-1$.
\end{main}








\subsection*{Computation of the Fourier coefficients}
Recall that Siegel's computation of the coefficients $a_h=a_h(E_k^{m})$ : 
$$
E_k^{m}(z)=\sum_{\gamma\in P\backslash \Gamma}\det(cz+d)^{-k}
=\sum_{h\in B_m}a_hq^h
$$
is based on the Poisson summation formula giving the equality
(see \cite{Maa71}, p.304): 
$$
\sum_{a\in S_m}\det(z+a)^{-k}={(-2\pi i)^{mk}\over 2^{m(m-1)\over 2}\Gamma_m(k)} \sum_{h\in C_m}\det(h)^{k-\frac {m+1} 2}
e^{2\pi i\tr (hz)},
$$
where $\Gamma_m(k)=\pi^{m(m-1)/4}\prod_{j=0}^{m-1}\Gamma(s-\frac j 2)$,
 $q^h=e^{2\pi i\tr (hz)}$, 
  $h$ runs over the set $C_m$ of positive definite half integral $m\times m$ symmetric matrices, and 
  $a$ runs over the set $S_m$ of integral $m\times m$ symmetric matrices,  see \cite{Si39}, p.652, \cite{St81}, p.338.
\subsection*{Formulas for the Fourier coefficients for $\det(2h)\not = 0$}  
\begin{align*}
a_h({E}_k^{m})&= 
{
(-2\pi i)^{mk}\Gamma_m^{-1}(k)
\over 
\zeta(k)\prod_{i=1}^{[m/2]}\zeta(2k-2i)
}
\\ & 
\times \det(2h)^{k-\frac {m+1} 2}M_h(k)
\begin{cases} \displaystyle 
L(k-\frac m 2,\psi_h ),
 & m \mbox{ even}, \\ 
1,& m \mbox{ odd}.  
\end{cases}
\end{align*}
The {\it integral factor}
  $ \displaystyle   M_h(k) = \prod_{{\ell} \in P(h)} M_{\ell}(h,{\ell}^{-k})$ 
  is a {\it  finite Euler product}, extended over primes ${\ell}$ in the set $P(h)$ of
  prime divisors 
  of all elementary divisors of the
  matrix $h$. The important property of the product is that for each ${\ell}$ we
  have that $M_{\ell}(h,t) \in \Z[t]$ is a polynomial with integral coefficients.
  
  Notice the {\it  $L$-factor $L(k-\frac m 2,\psi_h )$} depends on the  index $h$
of the Fourier coefficient;
this makes a difference to the case of odd $m$; the case of $GL(2)$ corresponds to $m=1$.
\subsection*{Proof: \small the use of the normalized Siegel-Eisenstein series}
defined as in \cite{Ike01}, \cite{PaSE} and [PaLNM1990] by \\ 
$$
{\mathcal E}_k^{m} =
{
{E}_k^{m}(z)
2^{m/2}
\zeta(1-k)\prod_{i=1}^{[m/2]}\zeta(1-2k+2i), 
}
$$
I show that it produces a nice $p$-adic family, namely:
\begin{prop}
\label{NSES}   \ 
\begin{itemize}
\item[(1)] 
For any non-degenerate matrice $h\in C_m$ the following equality holds
\begin{align}\label{ah}&
a_h({\mathcal E}_k^{m})=
2^{-\frac{m}2}\det h^{k-\frac {m+1} 2} M_h(k)
\\ & \nonumber
\times
\begin{cases} L(1-k + {\frac m 2}, \psi_h )C_h^{\frac m 2 - k + (1/2)},& m \mbox{ even}, \\ 
1,& m \mbox{ odd}, 
\end{cases} 
\end{align}
where $C_h$ is the conductor of $\psi_h$.
\item[(2)]  
for any prime $p>2$, and $\det(2h)$ not divisible by $p$, 
define the {\it $p$-regular part} $a_h({\mathcal E}_k^{m})^{(p)}$  of the coefficient 
$a_h({\mathcal E}_k^{{m}})$ of  ${\mathcal E}_k^{{m}}$  by introducing the factor
$\begin{cases} 
(1-\psi_h(p)p^{k-\frac m 2-1}){C_h^{\frac m 2 - k + (1/2)}},& m \mbox{ even}, \\ 
1,& m \mbox{ odd.}  
\end{cases}
$\\
Then $a_h({\mathcal E}_k^{m})^{(p)}$ 
is a $p$-adic analytic Iwasawa function of $t=(1+p)^k-1$ for all $k$ with $\omega^k$ fixed,
 and divided by the elementary factor $1-\psi_h(c_h)c_h^{k-\frac m 2}$. 
\end{itemize}
\end{prop}
\subsection*{Proof of (1) of Proposition \ref{NSES}}
Proof of (1) is deduced like at p.653 of \cite{Ike01}
from the Gauss duplication formula 
$$\Gamma (\frac s 2 )\Ga ( \frac {s+1} 2)=2^{1-s}\sqrt{\pi}\Gamma (s),$$
the definition 
$$
\Gamma_m(k)=\pi^{m(m-1)/4}\prod_{j=0}^{m-1}\Gamma(s-\frac j 2)
$$
and the functional equations
\begin{align*}&
\zeta(1-k)={2(k-1)!\over (-2\pi i)^k} \zeta(k), \\ &
\zeta(1-2k+2i)={2 (2k-2i-1)!\over (-2\pi i)^{2k-2i}} \zeta(2k-2i),
\\ & 
L(1-k+\frac m 2, \psi_h)
 ={2 (k-\frac m 2-1)!\over (-2\pi i)^{k-\frac m 2}} 
 L(k-\frac m 2, \psi_h) C_h^{k-\frac m 2-\frac 1 2} 
\end{align*}

\subsection*{Proof of (2) of Proposition \ref{NSES}} is then deduced easily : 

Notice that for any $a\in \Z_p^*$, the function of $t=(1+p)^k-1$
\begin{align}\label{ak}
k\mapsto a^k  &=  \omega(a)^k\langle a\rangle^k=
 \omega(a)^k 
 (1+p)^{k{\log \langle a\rangle\over \log(1+p)}}
 \\ \nonumber
 &
 =\omega(a)^k 
 (((1+p)^k-1)+1)^{\log \langle a\rangle\over \log(1+p)}
 \\& \nonumber
 =\omega(a)^k\sum_{n=0}^\infty{{\log \langle a\rangle \over \log(1+p)}\choose n}t^n  
\end{align} 
is a $p$-adic analytic Iwasawa function  denoted by $\tilde a(t)\in \Z_p[\![t]\!]$, of $t=(1+p)^k-1$ with $\omega^k$ fixed, where
${x\choose n}$=${x(x-1)\cdots (x-n+1)\over n!}$. 

Then
Mazur's formula applied to $L(1-k+\frac m 2,\psi_h)(1-\psi_h(p)p^{k-\frac m 2-1})$ shows that this function  is a $p$-adic analytic Iwasawa function of $t=(1+p)^k-1$ with $\omega^k$ fixed
(a single simple pole may occur at $k=\frac m 2$ only if $\omega^{k-\frac m 2}$ is trivial).

\

\subsection*{Proof of Main Theorem \ref{pMSE}}
Let us use the equality
$$
{E}_k^{m} =
{
\mathcal E}_k^{m}(z)\cdot
{2^{-m/2}\over
\zeta(1-k)\prod_{i=1}^{[m/2]}\zeta(1-2k+2i) 
}
$$
and the properties of the normalized series
${\mathcal E}_k^{n}(z)$ in Proposition \ref{NSES}.

First let us compute the reciprocal 
of the product of $L$-functions
$$
\zeta(1-k)\prod_{i=1}^{[m/2]}\zeta(1-2k+2i)
$$ 
using the above:
for even $k\ge 2$,  
\begin{align}\label{zk}&\zeta(1-k)^{-1}(1-p^{k-1})^{-1}  
={
U^*_{\theta_{k},c}((1+p)^{k-1}-1) (1-c^k)
\over P_{\theta_k,c}((1+p)^{k-1}-1)}\\ & \nonumber
\zeta(1-2k+2i)^{-1}(1-p^{2k-2i-1})^{-1}\\ &  \label{z2k}
={
U^*_{\theta_{2k-2i},c}((1+p)^{2k-2i-1}-1) (1-c^{2k-2i})
\over P_{\theta_{2k-2i},c}((1+p)^{2k-2i-1}-1)}
\end{align} 
which is meromorphic in the unit disc
with a finite number of poles (expressed via roots of $P_\theta$) for $\theta_k=\omega^{k-1}$.

\

Let us use again the notation $1+t=(1+p)^k$ with $k\in \Z_p$.
\begin{align}\label{pNF}& \nonumber
{2^{-m/2}\over
{\zeta(1-k)(1-p^{k-1})\prod_{i=1}^{[m/2]}\zeta(1-2k+2i)(1-p^{2k-2i-1})} 
}
\\ &
={{
U^E_{\omega^k}(t)
}
\over {
P^E_{\omega^k}(t)}}
\end{align}
where  the {\it numerator} is (an Iwasawa function)
$U^E_{\omega^k}(t)=$
$$U^*_{\theta_{k},c}(\frac{1+t}{1+p}-1) (1-c^k)
\prod_{i=1}^{[m/2]}
U^*_{\theta_{2k-2i},c}({\frac{(1+t)^2}{(1+p)^{2i+1}}-1}) (1-c^{2k-2i}), 
$$
and
$$
P^E_{\omega^k}(t)=
P_{\theta,c}\left({\frac{1+t}{1+p}}-1\right)
\prod_{i=1}^{[m/2]} P_{\theta_{2k-2i},c}\left({\frac{(1+t)^2}{(1+p)^{2i+1}}-1}\right)
$$
is the {\it polynomial denominator}  which depends only on $k \bmod p-1$.
\subsection*{Proof of Main Theorem \ref{pMSE}: control over the conductor of $\psi_h$}
Moreover, Mazur's formula applied to $L(1-k+\frac m 2,\psi_h)(1-\psi_h(p)p^{k-\frac m 2-1})$ ({\it in the numerator}) shows that
for all $h$ with $\det(2h)$ not divisible by $p$, 
\begin{align}\label{Lpsi}
L(1-k +\frac m 2, \psi_h)&(1-\psi_h(p)^{k-\frac m 2-1})\\ &\nonumber
={ 
G_{\theta, h}((1+p)^{k-\frac m 2-1}-1) 
\over
1-\psi_h(c_h)c_h^{k-\frac m 2}} 
\end{align} which is meromorphic in the unit disc 
with a possible single simple pole at $k=\frac m 2$  for all $k$ with $\theta=\omega^{k-1}$.
It comes from Mazur's measure on the finite  product $\displaystyle\prod_{{\ell}\in P_h}\Z_{\ell}^*$ extended over primes ${\ell}$ in the set $P_h=P(h)\cup \{p\}$; recall that $P(h)$ is the set of
  prime divisors 
  of all elementary divisors of the
  matrix $h$ as above. 

Indeed, for any choice of a natural number $c_h > 1$ coprime to $\prod_{{\ell}\in P_h}{\ell}$, there exists a $p$-adic measure 
$\mu_{c_h,h}$ on $\Z_p^*$, such that the special values
\begin{align*}&
L(1-k+\frac m 2, \psi_h)(1-\psi_h(p)p^{k-1-\frac m 2})=
{\int_{\Z_p^*}y^{k-\frac m 2-1}d\mu_{c_h, h}\over 
1-\psi_h(c_h)c_h^{k-\frac m 2}}
\\ & \nonumber
: =
(1-\psi_h(c_h)c_h^{k-\frac m 2})^{-1}
\int_{
\prod_{{\ell}\in P_h}\Z_{\ell}^*}\psi_h(y) y_p^{k-\frac m 2-1}d\mu_{c_h},  
\end{align*}
where Mazur's measure $\mu_{c_h}$ extends on the product 
$\displaystyle \prod_{{\ell}\in P_h}\Z_{\ell}^* {\buildrel y_p   \over \longrightarrow} \Z_p^*$ 
(see \S 3,Ch.XIII of \cite{LangMF}): 
$$
\mu_{c_h}(a+(N))
=
\frac {1}{c_h}\left[ \frac {c_ha}{N}\right] + 
\frac {1-c_h}{2c_h}=\frac{1}{c_h} B_1(\{\frac {{c_h}a} {N}\})-B_1(\frac {a}{N})
$$
for any natural number $N$ with all prime divisors in $P_h$.

 The regularizing factor is the following Iwasawa function which depends on $c_h \bmod 4p$ and $k\bmod p-1$:
\begin{align}\label{ck1}
&1-\psi_h(c_h) c_h^{k-\frac m 2} =1-(\psi_h\omega^{k-\frac m 2})((\frac{(1+p)^{k}}
{(1+p)^{\frac m 2}}
-1)+1)^{\log \langle c_h\rangle\over \log(1+p)}
 \\& \nonumber
 =1-(\psi_h\omega^{k-\frac m 2})(c_h)\sum_{n=0}^\infty{{\log \langle c_h\rangle \over \log(1+p)}\choose n}(
 (\frac{(1+p)^{k}}
{(1+p)^{\frac m 2}}
-1
 )^n
 \\& \nonumber
 =1-(\psi_h\omega^{k-\frac m 2})(c_h)\sum_{n=0}^\infty{{\log \langle c_h\rangle \over \log(1+p)}\choose n} (\frac{(1+t)}
{(1+p)^{\frac m 2}}
-1)^n
 \in \Z_p[\![t]\!],
\end{align}
where we write $c_h$ in place of $i_p(c_h)$ and use the notation $1+t=(1+p)^{k}$. 
The function (\ref{ck1}) is divisible by $t$ or invertible in $\Z_p[\![t]\!]$
according as $\omega^{k-\frac m 2}\psi_{c_h}$  is trivial or not
because $t=0 \Eq k=0$ and $1+t=(1+p)^k$.
\\

\subsection*{Elementary factors}

Notation: 
\begin{align*}
u_{c_h}(t)=
\begin{cases}
(1-\psi_h(c_h)\tilde c_h(t))/t, &\mbox{ if } \omega^{k-\frac m 2}\psi_h  \mbox{ is trivial}, \\
1- \psi_h(c_h) \tilde c_h(t),&  \mbox{ otherwise}.
\end{cases}
\end{align*}
By (\ref{ck1}) we have that  $u_{c_h}(t)\in \Z_p[\![t]\!]^*$, and we denote by $u^*_{c_h}(t)$ its inverse. 
Moreover,  (\ref{ah}) gives
the elementary factor 
$$
{\mathcal M}_h((1+p)^k-1)= 2^{-\frac{m}2}\det h^{k-\frac {m+1} 2} \prod_{{\ell}|P(h)} M(h, {\ell}^{-k})C_h^{k-\frac {m+1} 2}
$$
which is also an {\it Iwasawa function} as above:
 $$
 {\mathcal M}_h((1+p)^k-1)={\mathcal M}_h(t)\in \Z_p[\![t]\!].$$

\subsection*{Proof of Main Theorem \ref{pMSE}: the numerator}
It follows that
$$
a_h^{(p)}(k)=
{S^E_{\theta,h}((1+p)^{k}-1)\over P^E_{\theta,h}((1+p)^{k}-1)}={S^E_{\theta,h}(t)\over P^E_{\theta,h}(t)}, 
$$
where 
\begin{align*}&
S^E_{\theta,h}= u^*_{c_h}((1+p)^{k}-1)
{\mathcal M}((1+p)^{k}-1)
\\ & \times
 U_{\theta, h}((1+p)^{k-\frac m 2-1}-1) (1-c_h^k)
 U^*_{\theta_{k},c}((1+p)^{k-1}-1) 
 \\ & \times
\prod_{i=1}^{[m/2]}
U^*_{\theta_{2k-2i},c}((1+p)^{2k-2i-1}-1) (1-c_h^{2k-2i})
 \\ &
= u^*_{c_h}(t)
{\mathcal M}(t) U_{\theta, h}((1+t)(1+p)^{-\frac m 2-1}-1) 
\\ & \times
(1-\tilde c_h (t))
U^*_{\theta_{k},c_h}((1+t)(1+p)^{-1}-1) 
\\ & \nonumber 
\times \prod_{i=1}^{[m/2]}
U^*_{\theta_{2k-2i},c_h}((1+t)^2(1+p)^{-2i-1}-1) 
(1-\tilde c_h^{2}(t)c_h^{-2i}), 
\end{align*}
\subsection*{Proof of Main Theorem  (end)}
{\large \it The denominator is the following distinguished polynomial}

\begin{align*}&
P^E_{\theta, h} ((1+p)^{k-1}-1)= 
(1+((1+p)^{k-1}-2)\delta(\omega^{k-\frac m 2}\psi_{c_h}))
\\ & \times
P_{\theta_k,c_h}((1+p)^{k-1}-1) 
\prod_{i=1}^{[m/2]}
 P_{\theta_{2k-2i},c_h}((1+p)^{2k-2i-1}-1)
 \\ &
 = 
 (1+(t-1)\delta(\omega^{k-\frac m 2}
\psi_{c_h}))
P_{\theta_k,c_h}((1+t)(1+p)^{-1}-1) \\ & \nonumber
\times \prod_{i=1}^{[m/2]}
 P_{\theta_{2k-2i},c}((1+t)^2(1+p)^{-2i-1}-1), \mbox{ where}
\end{align*} 
$
\delta(\omega^{k-\frac m 2}\psi_{c_h})=
\begin{cases}
1, &\mbox{ if } \omega^{k-\frac m 2}\psi_{c_h}  \mbox{ is trivial}, \\
0,&   \mbox{ otherwise},  
\end{cases}
$ 
so that\\ 
$
1+(t-1)\delta(\omega^{k-\frac m 2}
\psi_{c_h})=
\begin{cases}
t, &\mbox{ if } \omega^{k-\frac m 2}\psi_{c_h}  \mbox{ is trivial}, \\
1,&   \mbox{ otherwise}.
\end{cases}
$ .\\
{\it
It remains to notice that different choices of $c_h$ coprime to $p \det (2h)$  give the same polynomial factors  $P^E_{\theta, h}$
 (up to invertible Iwasawa function). Indeed they all give the same single simple zero.
 }
\qed  

\

\section{Pseudomeasures and their Mellin transform }
\subsection*{Interpretation: Mellin transform of a pseudomeasure}
Pseudomeasures were introduced by J.Coates \cite{Co} as elements of the {\it fraction field ${\mathcal L}$}  of the 
Iwasawa algebra. 
Such a pseudomeasure is defined  by its Mellin transform 
which is a ring homomorphism and we can extend it by universality
(the extension of the integral along measures in $\Lambda=\Z_p[[T]]$ to the whole fraction field ${\mathcal L}$). 

The $p$-adic meromorphic function 
$$
a_h^{(p)}(k)=
{S^E_{\theta,h}((1+p)^{k}-1)\over P^E_\theta((1+p)^{k}-1)}={S^E_{\theta,h}(t)\over P^E_\theta(t)}, 
$$
is attached to an explicit pseudo-measure: 
$$
\rho_h^E={\mu_h^E \over \nu_h^E},\ \ 
 {S^E_{\theta,h}(t)\over P^E_\theta(t)}={\int_{\Z_p^*}\theta \chi_{(t)} \mu_h \over 
 \int_{\Z_p^*}\theta \chi_{(t)} \nu_h }
$$
\begin{itemize}
\item ${\mathcal S}(x)= \int_{\Z_p^*} x \mu_h^E$ is given by the collection of Iwasawa 
 functions $S_\theta(t)=\int_{\Z_p^*}\theta\chi_{(t)} \mu_h^E$ (the {\it numerator}),
\item  ${\mathcal P}(x)= \int_{\Z_p^*} x \nu_h^E$ is given by the collection of polynomials 
$P_\theta(t)=\int_{\Z_p^*}\theta\chi_{(t)} \nu_h^E$ (the {\it denominator}). 
\end{itemize}
   
\subsection*{Pseudomeasure $\rho$ as a family of distributions}
A pseudomeasure $\rho$ can be described as a certain {\it family of distributions, 
parametrized by the   set $X_p$ of $p$-adic characters}.

For any $x\in X_p$ 
we have a distribution given by the formula
$$
\rho_{h,x}^E (a+(p^v))=
\frac 1 {\varphi(p^v)} 
\sum^{}_{\hskip-.1cm\chi\bmod p^v}{\hskip-.3cm{}'}
\chi(a)^{-1}\frac {{\mathcal S}^E(\chi x)}
{{\mathcal P}^E(\chi x)}
$$
where ${}'$ means that the terms with ${\mathcal P}(\chi x)=0$ are omited.
It follows that
$$
\int_{\Z_p^*} \chi  \rho_{h,x}^E  
= \begin{cases}
\displaystyle 
{{\mathcal S}^E(\chi x) \over {\mathcal P}^E(\chi x)}
, &\mbox{ if } {\mathcal P}^E(\chi x) \not =0 \\
0, &\mbox{otherwise}, 
\end{cases}
$$
where 
\begin{align*}&
{\mathcal S}^E_h(\chi x)=S^E_{(\chi x)_\Delta, h}(\chi x(1+p)-1)
=S^E_{\theta,h}((1+p)^t-1), \\ &\theta=(\chi x)_\Delta, (\chi x)(1+p)=1+t, 
\\ &
{\mathcal P}^E_h(\chi x)
=P^E_{(\chi x)_\Delta, h}(\chi x(1+p)-1)
=P^E_{\theta,h}((1+p)^t-1). 
\end{align*}
\subsection*{A geometric construction: Siegel's method and duality 
 }
 For any Dirichlet character $\chi\bmod p^v$ 
 consider Shimura's "involuted" Siegel-Eisenstein series
 assuming their absolute convergence (i.e. $k> m+1$):
 $$
E^*_k(\chi,z)=\sum_{(c,d)/\sim}\chi(\det(c))\det(cz+d)^{-k}
=\sum_{h\in B_m}a_h(E^*_k(\chi,z))q^h
$$
The series on the left is {\it geometrically defined}, and the Fourier coefficients on the right can be computed by Siegel's method (see \cite{St81}
\cite{Shi95}, \dots )
via the singular series
\begin{align}\label{ahk3}
&a_h({E}^*_k(\chi,z))\\ &\nonumber = 
{
(-2\pi i)^{mk}\over 2^{m(m-1)\over 2}\Gamma_m(k)} 
\sum_{{\mathfrak R}\bmod 1}\chi(\nu({\mathfrak R}))\nu({\mathfrak R})^{-k}
\det h^{k-\frac{m+1} 2}e_m(h{\mathfrak R})
\end{align}


If $\chi=\chi_0\bmod p$ is trivial and $p\nd \det h$ then
\begin{align}\label{ahk4}
&a_h({E}^*_k(\chi_0,z))\\ &\nonumber = 
{
(-2\pi i)^{mk}\over 2^{m(m-1)\over 2}\Gamma_m(k)} 
\sum_{{\mathfrak R}\bmod 1}\chi_0(\nu({\mathfrak R}))\nu({\mathfrak R})^{-k}
\det h^{k-\frac{m+1} 2}e_m(h{\mathfrak R})=
\\ & \nonumber
a_h(E_k^m)\times
\begin{cases} (1-p^{-k})(1+\psi_h (p)p^{-k+\frac m 2})\prod_{i=1}^{(m/2)-1}(1-p^{-2k+2i}),
& m \mbox{ even}
\\ \\
(1-p^{-k})\prod_{i=1}^{(m-1)/2}(1-p^{-2k+2i}),& m \mbox{ odd}.  
\end{cases}
\end{align}  

The formula (\ref{ahk4}) means that the series
${E}^*_k(\chi_0,z)$  
coincides with $E_k^m$ after
removing $h$ with $\det h$ divisible by $p$ and normalizing by the factor in
(\ref{ahk4}).
Moreover, the {\it Gauss reciprocity law} shows that
the normalizing factor depends only on $\det h\bmod 4p= \det h_0\bmod 4p$, where 
$h_0\equiv h \bmod 4p$ runs through a representative system.
Let us denote this factor by $C^+(h_0,k,4p)$: 
for the trivial character $\chi=\chi_0\bmod p$ and $\det h$ not divisible by $p$
\begin{align}\label{ahk5}
&a_h({E}^*_k(\chi_0,z))= 
a_h(E_k^m)C^+(h_0,k,4p),  \mbox{ where }\\ & \nonumber
C^+(h_0,k,4p)=
\begin{cases} (1-p^{-k})(1+\psi_h (p)p^{-k+\frac m 2})\prod_{i=1}^{(m/2)-1}(1-p^{-2k+2i}),
& m \mbox{ even}
\\ \\
(1-p^{-k})\prod_{i=1}^{(m-1)/2}(1-p^{-2k+2i}),& m \mbox{ odd}.  
\end{cases}
\end{align}  
{\large \it From the Fourier coefficients to modular forms: }

If we remove in the Fourier expansion 
$\it
E_k^m(z)=\sum_{h\ge 0}a_he_m(hz)
$ 
all terms with $\det h$ divisible by $p$ the equality of Fourier coefficients (\ref{ahk5}) transforms to the equality of the series

\begin{align}\label{chi0}&
{E}^*_k(\chi_0,z)=
(4p)^{-m(m+1)/2}\sum_{h_0\bmod 4p\atop p\nds \det h_0}C^+(h_0,k,4p)\times
\\ & \nonumber
\sum_{x\in S\bmod 4p}
e_m(-h_0x/4p)E_k^m(z+(x/4p)).
\end{align}
\subsection*{A geometric construction 
 }
Let us apply the interpolation theorem  (Theorem \ref{pMSE})
to all the coefficients 
\begin{align}\label{chi02}
& a_h^{(p)}(k)=a_h(E_k^{m})C^-(h_0,k,4p), \mbox{ where }
\\ & \nonumber C^-(h_0,k,4p)=
\begin{cases} 
{1-\psi_h (p)p^{k-\frac m 2-1}\over(1-p^{k-1})\prod_{i=1}^{m/2}(1-p^{2k-2i-1})},& m \mbox{ even}
\\  \\
{1\over(1-p^{k-1})\prod_{i=1}^{(m-1)/2}(1-p^{2k-2i-1})},& m \mbox{ odd},  
\end{cases}
\end{align}
and (\ref{chi0}) becomes a "geometric-algebraic equality"
of two families of modular forms 
\begin{align}\label{chi01}&
{E}^*_k(\chi_0,z)=
(4p)^{-m(m+1)/2}\sum_{h_0\bmod 4p\atop p \nds \det h_0}
C^+(h_0,k,4p)
\times
\\ & \nonumber
C^-(h_0,k,4p)\sum_{x\in S\bmod 4p}
e_m(-h_0x/4p)E_k^m(z+(x/4p)).
\end{align}
\subsection*{A geometric construction (end)}
We deduce by the orthogonality that 
\begin{align}\label{chi03}&
\sum_{x'\in S\bmod 4p}
e_m(-h_0x'/4p)E^*_k(\chi_0, z+(x'/4p))=
\\ \nonumber &
C^+(h_0,k,4p)
\sum_{x\in S\bmod 4p}
e_m(-h_0x/4p)
E_k^m(z+(x/4p)).
\end{align}
Each series
$
C^-(h_0,k,4p)\sum_{x\in S\bmod 4p}
e_m(-h_0x/4p)E_k^m(z+(x/4p))
$
in (\ref{chi01})
 determines a {\it unique pseudomeasure with coefficients in 
${\bar\Q}[\![q^{B_m}]\!]$} whose moments are given by those of the coefficients (\ref{chi02}).The unicity means that a pseudomeasure is determined by its Mellin transform.
It is also a family of distributions geometrically defined by the series
$$
{C^-(h_0,k,4p)\over C^+(h_0,k,4p)}
\sum_{x\in S\bmod 4p}
e_m(-h_0x/4p)E^*_k(\chi_0, z+(x/4p)).
$$ 

\


\section{Application to Minkowski-Siegel Mass constants}
\subsection*{$p$-adic version of Minkowski-Siegel Mass constants.
}
An application of the construction is the  $p$-adic version of Siegel's Mass formula.
It expresses the Mass constant through the  above product of $L$-values.
This    
product can be viewed as the proportionality coefficient between two kinds of Eisenstein series in the symplectic case extending Hecke's
result (1927) 
of the two kinds of Eisenstein series and 
the relation between them. 
However, there is no direct analogue of Hecke's computation in the symplectic case.

Thus this mass constant admits an 
explicit product expression through the values of the functions (\ref{zki}) at $t_j=(1+p)^j-1$, 
for $j=k$, and $j=2,4,\dots, 2k-2$.


\

Recall that (\cite{ConSl98}, p.409)
\begin{quote}
unimodular latticies have the property that there are explicit formulae, the mass formulae, which give appropriately weighted sums of the theta-series of all the inequivalent latticies of a given dimension.
In particular,  the numbers of inequivalent latticies is given by Minkowski-Siegel Mass constants for unimodular latticies.
\end{quote} 

\



In the particular case of even unimodular
quadratic forms of rank $m = 2k \equiv 0(\bmod 8)$, this formula means that
there are only finitely many such forms up to equivalence for each $k$ and that,
if we number them $Q_1, \dots ,Q_{h_k}$ , then we have the relation
$$
\sum_{i=1}^{h_k}\frac 1 {w_i}\Theta_{Q_i}(z)=m_kE_k
$$
where $w_i$ is the number of automorphisms of the form $\Theta_{Q_i}$ is the theta
series of $Q_i$, $E_k$ the normalized Eisenstein series of weight $k = m/2$ (with
the constant term equal to 1),

The dimension of lattices is $2k$ and the Mass formula 
express  an identity of a sum of {\it weighted theta functions}
 and a {Siegel-Eisenstein series} of weight $k$, 
 multiplied by the Mass constant 
 $$
 m_k=2^{-k}
\zeta(1-k)\prod_{i=1}^{k-1}\zeta(1-2k+2i)= 
(-1)^{k}\frac {B_k} {2k} 
\times \prod_{j=1}^{k-1}\frac {B_{2j}} {4j}
 $$
which is related 
the above normalising coefficient.\\
{\tt 
gp > mass(4)\\
\% = 1/696729600\\
gp > mass(8)\\
\% = 691/277667181515243520000
}

\


The present result says that the $p$-regular part of $1/m_k$ is a product of values of the 
$p$-adic meromorphic functions (\ref{zki}) at $t_j=(1+p)^j-1$, $j=k$ and $j=2, 4, \dots, 2k-2$.  

It is known that the rational number $m_k$ becomes very large rapidly, 
when $k$ grows (using the functional equation). 
It means that the denominator of $1/m_k$ becomes {\it enormous}.
The explicit formula (\ref{pNF}) 
applied to  the reciprocal of the product of $L$-functions as above 
shows that these are {\it only irregular primes} which contribute to the denominator, 
and this contribution can be evaluated for {\it all primes} knowing the {\it Newton polygons} of the polynomial part $P_\theta$, which can be found directly from 
the Eisenstein measure.
Precisely, for the distingushed polynomial $P(t)=P_{\theta}(t)=a_dt^d+\dots+a_0$, 
$\ord_p a_d=0$, and $\ord_p a_i>0$ for $0\le i\le d-1$, and $\ord_p(t_j)=\ord_p j+1$, where $t_j=(1+p)^j-1$  for $j=k$ and $j=2, 4,\dots,2k-2$.
Then 
$$
\ord_p P(t_j)= \mathop{\min}_{i=0, \dots, d} \left(\ord_p a_{i,k}+i(\ord_p j+1)\right).
$$
the values $\ord_p a_{i,k}$ for $0\le i\le d$ come from the Iwasawa series in the denominator in the left hand side of (\ref{pNF}).
Also, it gives an important information about the location
 of zeroes  of the polynomial part as in (\ref{pNF})). 
However $P(t_j)\ne 0$ in our case because all the $L$-values in question do not vanish. 

 \

\subsection*{Application to Minkowski-Siegel Mass constant \\ (numerical illustration)}
{\tt
for(k=1,10,print(2*k,~factor(denominator(1/mass(2*k)))))\\
{2}       1\\
4       1\\
6       1\\
8       [691, 1]\\
10       [691, 1; 3617, 1; 43867, 1]\\
12       [131, 1; 283, 1; 593, 1; 617, 1; 691, 2; 3617, 1; 43867, 1]\\
14       [103, 1; 131, 1; 283, 1; 593, 1; 617, 1; 691, 1; 3617, 1; 43867, 1; 6579
31, 1; 2294797, 1]\\
16       [103, 1; 131, 1; 283, 1; 593, 1; 617, 1; 691, 1; 1721, 1; 3617, 2; 9349,
 1; 43867, 1; 362903, 1; 657931, 1; 2294797, 1; 1001259881, 1]\\
18       [37, 1; 103, 1; 131, 1; 283, 1; 593, 1; 617, 1; 683, 1; 691, 1; 1721, 1;
 3617, 1; 9349, 1; 43867, 2; 362903, 1; 657931, 1; 2294797, 1; 305065927, 1; 
 1001259881, 1; 151628697551, 1]\\
20      [103, 1; 131, 1; 283, 2; 593, 1; 617, 2; 683, 1; 691, 1; 1721, 1; 3617,
1; 9349, 1; 43867, 1; 362903, 1; 657931, 1; 2294797, 1; 305065927, 1; 
1001259881, 1; 151628697551, 1; 154210205991661, 1; 26315271553053477373, 1]
}

\

\section{Link to Shahidi's method for 
 SL(2)  and regular prime $p$}
\subsection*{Methods of  constructing $p$-adic $L$-functions}
{Our long term purposes are }
to define and to use the $p$-adic  $L$-functions in a way similar to complex  $L$-functions
via the following methods:

(1) Tate, Godement-Jacquet;

(2) the method of Rankin-Selberg;

(3) the method of Euler subgroups of  Piatetski-Shapiro and the  doubling method of  Rallis-Böcherer (integral representations  on a subgroup of $G\times G$);

(4) Shimura's method (the convolution integral with theta series), and

(5) Shahidi's method.

There exist already advances for (1) to (4), and we are also trying to develop  (5).

We use the Eisenstein series on classical groups and  $p$-adic integral of Shahidi's type
 for the reciprocal of a product of certain $L$-functions.
 
 \
 

\subsection*{Link to Shahidi's method in the case of 
 SL(2)  and regular prime $p$}
 
The starting point here is the Eisenstein series
$$
E(s,P,f,g)=\sum_{\gamma\in P\smallsetminus G}
f_s(\gamma g),
$$
on a reductive group $G$  and a {\it maximal  parabolic subgroup} $P=MU^P$
 (decomposition of Levi). 

This series generalizes
$$
E(z,s)=\frac{1}{2}\sum\frac{y^s}{|cz+d|^{2s}},\quad (c,d) =1.
$$

Here $f_s$ is an appropriate function in the induced representation space 
 $I(s, \pi)=Ind_{P_\A}^{G_\A}(\pi\otimes |det_{M}(\cdot)|^s_{\A}) )$, see
 (I.2.5.1) at p. 34 of \cite{GeSha}.
 
\subsection*{
Computing a non-constant term  (a Fourier coefficient)}
 of this Eisenstein series provides an analytic continuation and the functional equation  for many Langlands
 $L$ functions $L(s,\pi,r_j)$. 
 
 \
 
In this way the  $\psi$-th Fourier coefficient  (with $\psi$ of type $\psi(x)=\exp(2\pi i n x), n\in \N, n\ne 0$) of the series $E(s,P,f,e)$ is determined by the  Whittaker functions $W_v$  in the form (see \cite{GeSha}, (II.2.3.1), p.78):
$$
E^\psi(e,f,s)= \prod_{v\in S}W_v(e_v)\prod^m_{j=1}\frac{1}
{L^S(1+js,\pi, {r}_j)},
$$
where ${r}_j$ are certains fundamental representations of the dual group ${}^LM$.
\subsection*{}
\begin{theo}[a complex version] 
With the data $G=SL(2)$, $M =\{\mat a, 0,  0, a^{-1}, \}\cong \GL_1 $, $\pi=I$, and $\psi$ a  non-trivial character of the group 
$U(\ab)/U(\qb)$, $U=\{\mat 1, *,  0, 1, \}\cong \G_a$, 
let
$
E^{\psi}(s,f,e)= \int E(s,f,n){\psi(n)} dn,
$
the integration on the  quotient  space of $U(\ab)$ by $U(\qb)$.
Then the first Fourier coefficient has the form
$$
E^{\psi}(s,f,e)={W_{\infty}}(s)\frac{1}{\zeta(1+s)}, 
$$
for a certain Whittaker function ${W_{\infty}}(s)$ (see \cite{Kub}, p.46).
\end{theo}
\begin{theo}[a $p$-adic  version, 
a work in progress] (with S.Gelbart, S.Miller, 
F.Shahidi) 

Let $p$ be a regular prime. 
Then there exists an explicitly given distribution   $\mu^*$ on  $\zb^{*}_p$ 
such that for all   $k \geq 3$ and for all primitive Dirichlet characters $\chi\bmod p^v$ with $\chi(-1)=(-1)^k$ one has 
$$
\int_{\zb^{*}_p} \chi y^{k}_{p} \mu^*=\frac{1}{(1-\chi(p) p^{k-1})L(1-k, \chi)}, 
$$
where  $L(s, \chi)$ 
is the Dirichlet $L$-function.
 More precisely, the distribution $\mu^{*}$ can be expressed through the non-constant Fourier coefficients of a certain   Eisenstein series $\Phi^{*}$.
\end{theo}
\hskip-3pt{\large \it Remark.}   Using Siegel's method for the symplectic groups
$GSp_m$,  and for all primes $p$, this result also follows from Main Theorem \ref{pMSE}
by specializing it to the case of regular $p$ and $m=1$.


\

\section{
Doubling method and 
Ikeda's constructions}
\subsection*{Further applications: we only mention the proof of the $p$-adic Miyawaki Modularity Lifting Conjecture}
 by pullback of 
families  
Siegel modular forms (jointly with
Hisa-Aki Kawamura), see \cite{Kawa}, \cite{PaIsr11}.

\

Ikeda's constructions (\cite{Ike01}, \cite{Ike06}) extend the {doubling method} to pullbacks of cusp forms instead of pullbacks of Eisenstein series. 

\

In the Fall 1999 in IAS, Ilya was much inspired by the preprint of the first Ikeda's  lifting,  and tried to interpret it representation-theoretically.

Indeed, it extends his own work \cite{PS1} on Saito-Kurokawa lifting from genus 2 to arbitrary genus $2m$.

In fact, there is a relation of Ikeda's work to Arthur's conjecture \cite{Ar89}.  

\

In the same period, Ilya studied  the preprint of 
{\it \cite{KMS2000} on $p$-adic Rankin-Selberg $L$-functions in an informal  seminar in his office in IAS
together with me and other participants:  Jim Cogdell, Siegfried Böcherer, Reiner Schulze-Pillot, \dots} 

\

\subsection*{The use of the 
The Eisenstein family 
${\mathcal E}_k^{(n)}$} as above plays a crucial role in Ikeda's work:
the idea was to substitute the Satake parameter $\alpha_p(k)$ of a cusp form in place of the 
parameter $k$ in the Siegel-Eisenstein family.  

Both $p$-adic and complex analytic $L$-functions are produced in this way.

Thus obtained cuspidal $p$-adic measures generalize the Eisenstein measure, and produce families of cusp forms.

A version of this construction produces Klingen-Eisenstein series 
 and Langlands Eisenstein series, see \cite{PaSE} ($p$-adic Peterson product of a cusp form with a pullback of the constructed family), more recently used by Skinner-Urban \cite{MC}.

\

For genus two, my student P.Guerzhoy found  in 1998 a p-adic version of the {\it holomorphic Maass-Saito-Kurokawa lifting} \cite{Gue}, answering a question of E.Freitag.  
P.Guerzhoy visited Ilya here in Yale in 1999.

\def\mat #1,#2,#3,#4,{\left({#1\atop #3}{#2\atop #4}\right)}
\def\frac #1,#2,{{{#1}\over {#2}}}
\def\bra#1,{{\left\lbrace {#1}\right\rbrace}}
\def\ph{\varphi}
\def\Ph{\Phi}
\def\si{\sigma}
\def\z{\zeta}
\def \sign {{\rm sign}}
\def \sgn {{\rm sgn}}
\def \inv {{\rm inv}}
\def \diag {{\rm diag}}
\def \rem {{\rm rem}}
\def \Br {{\rm Br}}
\def \Inf {{\rm Inf}}
\def \Proj {{\rm Proj}}
\def \Mor {{\rm Mor}}
\def \Div {{\rm Div}}
\def \div {{\rm div}}
\def \Pic {{\rm Pic}}
\def \Res {{\rm Res}}
\def \vol {{\rm vol}}
\def \Ver {{\rm Ver}}
\def \di {{\rm div}}
\def \Spin {{\rm Spin}}
\def \Step {{\rm Step}}
\def \Spl {{\rm Spl}}
\def \M{{\rm M}}
\def \Gal{{\rm Gal}}
\def \suc{{\rm suc}}
\def \Gab{G^{\rm ab}}
\def \rk{{\rm rk}}
\def\tm{\tilde m}
\def\tGa{\tilde{\Ga}}

\def\ttl{\tilde t_{\la}}
\def\l1{\langle} 
\def \rr{\rangle}

\def\Ci{C^{\infty}}
\def\tr{{\rm tr}}
\def\Sh{{\ru X}\,}
\def\Tr{{\rm Tr}}
\def\Id{{\rm Id}}
\def\Re{{\rm Re}\,}
\def\Res{{\rm Res}}
\def\Im{{\rm Im}\,}
\def\Ind{{\rm Ind}}
\def\Irr{{\rm Irr}}
\def\Aut{{\rm Aut}}
\def\Ker{{\rm Ker}\,}
\def\End{{\rm End}}
\def\GSp{{\rm GSp}}
\def\Fr{{\rm Fr}}
\def\Char{{\rm Char}}
\def\Spec{{\rm Spec}}
\def\Card{{\rm Card}}
\def\Hom{{\rm Hom}}
\def\Aut{{\rm Aut}}
\def\pr{{\rm pr}\,}
\def\Or{{\rm O}}
\def\M{{\rm M}}
\def\qt{{\rm qt}}
\def\rad{{\rm rad}}
\def\Gd{{\rm Gd}}
\def\gd{{\rm gd}}

\def\nd{\not\hskip2.2pt\mid}
\def\th{{\theta}}
\def\D{{\cal D}}
\def\L{{\cal L}}
\def\Hc{{\cal H}}
\def\GS{{\rm GSp}_m}
\def\det{{\rm det}\,}
\def\Supp{{\rm Supp}}
\def\mod{{\rm mod}}
\def\H{{\mathfrak H}}
\def\Sp{{\rm Sp}}
\def\GL{{\rm GL}}
\def\G{{\bf G}}
\def \i{{\bf i}}

\def\Qa{\Q^{{\rm ab}}}
\def\Zp{({\Z}^+)}
\def\al{\alpha}
\def\w{\omega}
\def\W{\Omega}
\def\ga{\gamma}
\def\Ga{\Gamma} 
\def\la{{\lambda}}
\def\La{{\Lambda}}
\def\ka{{\kappa}}
\def\N{{\bf N}}
\def\No{{\rm N}}
\def\O{{\cal O}}
\def\SL{{\rm SL}}
\def\Ind{{\rm Ind}} 
\def\De{{\Delta}}
\def\de{{\delta}}
\def\Si{{\Sigma}}
\def\V{{\cal V}}
\def\ep{{\varepsilon}}
\def\is{\buildrel \sim   \over \rightarrow}
\def\mi{\buildrel \cdot   \over -}
\def\sqrt{\radical"270370 }
\def\Hol{{\cal H}ol\,}
\def\Distr{{\cal D}istr\,}
\def\oe{\overline }
\def\M{{\rm M}}
\def\Mc{{\cal M}}
\def\Sc{{\cal S}}
\def\Fc{{\cal F}}
\def\Ec{{\cal E}}
\def\Lc{{\cal L}}
\def\Dc{{\cal D}}
\def\Mp{M^\prime}
\def\Np{N^\prime}
\def\Ch{C_{\chi}}  
\def\Cc{{\cal C}}
\def\Nc{{\cal N}}
\def\Vc{{\cal V}}
\def\a{{\go a}}
\def\b{{\go b}}

\def \d{{\go d}}
\def \n{{\go n}}
\def \m{{\go m}}
\def \mp{{\go m}^\prime}
\def \p{{\go p}}
\def \q{{\go q}}
\def \e{{\go e}}

\def \f{{\bf f}}
\def \Ng{{\go N}}
\def \Pg{{\go P}}
\def \g{{\bf g}}
\def \P{{\bf P}}
\def \GL{{\rm GL}}
\def \SI{{\rm Spin}}
\def \SO{{\rm SO}}
\def \GS{{\rm GSp}} 
\def  \X{{\cal X}}
\def \h {{\go H}}
\def \tmu {{\tilde {\mu}}}
\def \GA{G_{\A}}
\def \GQ{G_{\Q}}
\def \Hc{{\cal H}}
\def \Xc{{\cal X}}
\def \Ac{{\cal A}}
\def \np{n_1^\prime}
\def \npp{n_1''}
\def \ord{{\rm ord}}
\def \F {{\bf F}}
\def \Tp {T^{\prime}}
\def \ap {\a^{\prime}}
\def \Eq {\Longleftrightarrow}
\def \bs {\backslash}
\def\leaderfill{\leaders\hbox to 1em{\hss.\hss}\hfill}

\def \Eq {\Longleftrightarrow}
\def \bs {\backslash}
\def\leaderfill{\leaders\hbox to 1em{\hss.\hss}\hfill}

\renewcommand{\theequation}{\thesubsection.\arabic{equation}}
\newcommand{\Subsection}{\setcounter{equation}{0} \subsection}
\newpage
\appendix
\Section
{Appendix. On  $p$-adic $L$-functions for $GSp(4)$}
talk by Alexei Panchishkin
on December 2, 1999, at Automorphic Forms and $L$-functions Seminar in IAS.


\setcounter{subsection}{-1}
\Subsection{Introduction.}
The purpose of this talk is to describe a joint work in progress with
 I.I.Piatetski-Shapiro started in February 1998 in Jerusalem during the
 conference "$p$-Adic Aspects of the Theory of Automorphic Representations".

\medskip\noindent
Let $G$ be a semi-simple algebraic group over a number field $F$,
 and $p\ge 5$ be a fixed prime number. 
Recall that the {\sl Iwasawa albebra} $\La$ is defined as $\Z_p[[T]]$  and let
 $\Lc = {\rm Quot}
 \La$ denote its quotient field.
Elements $a(T)\in \Lc$ represent some $\C_p$-meromorphic functions with finite
 number of poles
 on the unit disc $U_p=\{t\in \C_p\ |\ |t|_p<1\}\subset \C_p$ where
$\C_p=\hat{\overline \Q}_p$ the {\sl Tate field}. 
We consider the following problem: how to attach
 to a (complex valued) Langlands $L$-functions $L(s, \pi, r)$ a certain
 $p$-adic valued meromorphic $L$ function 
 $L_{\pi, r, p}$  with a finite number of poles
  where $\pi$ is an automorphic representation of the adelic group $G(\A_F)$ and 
 $r$ is a finite dimensional complex representation $r: {}^LG(\C)\to GL_m(\C)$
 of the Langlands group ${}^LG(\C)$.
The $p$-adic $L$-function $L_{\pi, r, p}$ should belong to  $\Lc$ or to its
 finite extension. 
The first example of a function of this type comes from the work of
 Kubota and Leopoldt [Ku-Le] and Iwasawa [Iw]:
 there exists a unique element $g(T)\in \Lc$ such that for all positive
 integers $k\equiv 0 (\bmod (p-1))$, $g((1+p)^k-1)=\z^*(1-k)$, where
 $\z^*(s)=(1-p^{-s})\z(s)$ is the Riemann zeta function with the $p$-factor
 removed from its Euler product.
The function  $\z_p(s)=g((1+p)^{1-s}-1)$ is
 analytic for all $s\in \Z_p\backslash 1$ with values in $\Q_p$ and it is called
 the {\sl Kubota-Leopoldt $p$-adic zeta function}. 
It has the following properties: $\z_p(1-k)=\z^*(1-k)$ for all positive integers
 $k\equiv 0 (\bmod (p-1))$, and $\Res_{s=1}\z_p(s)=1-\frac 1,p,$.
In this case we have actually $Tg(T)\in \La^\times$ so that $\z_p(s)$ has no
zeroes, unnike the complex zeta-function.
However, one could start from another progrssion $k\equiv i(\bmod p-1), k>0$, 
 $i\bmod (p-1)$ and obtain in the same way other branches $\z_{p,i}(s)$
 of $p$-adic zeta function which have interesting zeroes important in the
 Iwasawa theory [Iw, Wi90].

\medskip\noindent
Constructions of $L_{\pi, r}\in \Lc$ are known in a number of
 cases but there exists no general definition.
For example, the standard $L$ functions $L(s, \pi, St_{2n+1})$
 of degree $2n+1$ for the group $GSp_{2n}\subset GL_{2n}$ over $F=\Q$
 attached to the standard
 orthogonal representation of ${}^L\GSp_{2n}(\C)$ and to a cuspidal
 irreducible representation $\pi=\pi_f$ coming from a
 holomorphic Siegel cusp eigenform $f$ admits a $p$-adic
 analogue which was constructed using the {\sl Rankin-Selberg method} in the
 $p$-ordinary
 case [PaLNM] for {\sl even} $n$.
This construction was extended by S.B\"ocherer and C.-G. Schmidt [Bo-Sch]
 to the general case of $p$-ordinary forms of arbitrary genus $n$ and weight 
 $k>n$, by using the {\sl method of doubling of variables}.
The critical values in the sense of Deligne [De79] of the $L$-function
 $L(s, \pi_f\otimes \chi, St_{2n+1})$ are $s\in \Z$ such that $1-k+n\le s \le
 k-n$
 satisfying the following parity condition: $\displaystyle{(-1)^s=
 \begin{cases}
\chi(-1),&\mbox{ if }s\ge 0 \\ -\chi(-1),&\mbox{ if }s<0
 \end{cases}
 }
 $ for a Dirichlet character
$\chi\bmod p^N$. 
This description follow from the form of $\Gamma$-factor 
\begin{align}&
L_\infty (s, \pi_f\otimes \chi,
St_{2n+1})=\prod_{j=1}^{n}\Ga_\C(s+k-j)\Ga(s+\delta) ,\\ & \nonumber
\Ga_\C(s)=2(2\pi)^{-s}\Ga(s), \Ga_\R(s)=(\pi)^{-s/2}\Ga(s/2), \de=(1-\chi(-1))/2.
\end{align} 
In this case the algebraic
 numbers  $\displaystyle{L_{s,\chi}=\frac {L^*((s, \pi_f\otimes
\chi, St_{2n+1})},{\langle f, f\rangle},}$ can be interpolated to values of some
Iwasawa-type series $g_{f, i}(\chi (1+p)(1+p)^k-1)$ where $\langle f, f\rangle$
is the Petersson scalar product, $i$ runs over residues $\bmod (p-1)$.
In this case ${}^LGSp_{2n}(\C)=GSpin_{2n+1}$, the universal cover of the 
 orthogonal group $GO_{2n+1}(\C)$, $St_{2n+1}:GO_{2n+1}(\C)\hookrightarrow
 GL_{2n+1}(\C)$.

 \medskip\noindent
In order to construct in general $p$-adic automorphic $L$-functions 
 out of their complex critical special values one can successfully use $p$-adic 
 integration along a (many variable) Eisenstein measure which 
 was introduced by N.Katz [Ka78]
 and used by H.Hida [Hi91] in the case of $G=\GL_2$ over a
 totally real
 field $F$ (i.e. for the elliptic modular forms and
 Hilbert modular forms).  
The application of such a measure to a given $p$-adic family 
 of modular forms provides a general construction of
 $p$-adic $L$-functions of several variables. 
On the other hand, the evaluation 
 of this measure at certain points gives another important source of 
 $p$-adic $L$-functions [Ka78]. 
In the Siegel modular case the Eisenstein measure was constructed in [PaSE]. 

\medskip\noindent
The goal of our work is to construct a $p$-adic version of the
 $L$-function $L(s, \pi_f, r_4)$ of degree 4 attached to a Siegel-Hilbert
cusp egenform of degree 4 over a totally real field $F$, i.e. for the 
 symplectic group  
$$
\GSp_4 =\left\lbrace
 g\in \GL _{4}\left\vert\right. {}^t g J_4 g =\nu(g)J_4,
 \nu (g ) \in \GL_1
\right\rbrace,
$$  
 over $F$ where 
$$
J_4=\begin{pmatrix}0_2&-1_2\cr1_2&0_2\end{pmatrix}
$$
We use the Eisenstein measure and a $p$-adic analogue of the Petersson product
 for $\La$-adic automorphic forms on
 $GL_2$ over a totally real field, see [Hi90, Hi94].
Instead of $p$-adic interpolation of critical values we try to imitate in the
 $p$-adic case a known complex analytic integral representation for
 $L(s, \pi_f, r_4)$.  Main Theorem is given in Section 4.

\Subsection{
Complex analytic  $L$-functions for $GSp(4)$.}
Let $F$ be a global field of characteristic $\not = 2$, and $V$ a four
 dimensional vector space over $F$ endowed with  a non-degenerate
 skew-symmetric form $\rho : V\times V\to F$, 
$$
G_\rho=\GSp_4 =\left\lbrace
g\in \GL (V)\left\vert\right. \rho (gu, gv)=\nu_g\rho(u, v), \ 
\nu_g\in F^\times
\right\rbrace,
$$ 
 the algebraic group of symplectic similitudes of $\rho$ 
 over $F$.
Let $\pi=\otimes_v\pi_v$ be an irrreducible cuspidal automorphic representation
 of $G_\rho(\A_F)$ where $v$ runs over all places of $F$, then according to
 Langlands' classification of irreducible supercuspidal representations $\pi_v$
 of $G_\rho (F_v)$ for almost all $v$  $\pi_v$ correspond  to a semi-simple
 conjugacy class of a diagonal matrix 
\begin{align*}
h_v=\diag \{\al_0, \al_0\al_1, \al_0\al_2,
\al_0\al_1\al_2 \}\in &{}^LG_\rho (\C) \is GSP_4(\C) 
\buildrel {r_4} \over\to GL_4(\C)\cr &
(\al_j=\al_j(v), v\not \in S, |S|<\infty).
\end{align*}

The Andrianov $L$-function (or the {\sl spinor $L$-function}) of $\pi$ is then  
 the following Euler product 
\begin{align}
L(s, \pi, r_4)= \prod _{v\not \in S}\det (1_4-r_4(h_v)\cdot Nv^{-s})^{-1}\times
\left({{\rm a \ finite\  Euler\ product}}\atop {{\rm over \ }}v\in S \right)
\end{align}
This $L$ function plays an impotant role in arithmetic, in particular it is
 related  to $l$-adic Galois representation on $H^3$ of the
 corresponding Siegel threefold [Tay], [Lau].

\medskip\noindent
This $L$ function was introduced by Andrianov [AndBud], [And74]
 in  the classical fashion, for
 $F=\Q$, and for $\pi =\pi_f$ coming from a holomorphic Siegel cusp eigenform 
 $f=\sum_\xi A_\xi q^\xi$
 for the Siegel modular group $\Ga_2=Sp_4(\Z)$
 over the Siegel upper half plane of genus two
$$
H_2=\{ z={}^tz\in M_2(\C)\ |\  \Im (z)>0\},
$$
 where $\xi$ runs over the semi-group $B_2$ of semi-definite half integral 
 symmetric $2\times 2$-matrices $\xi$, $A_\xi\in \C$, so that
 $q^\xi=\exp(2 \pi i \Tr(\xi z))$ form a multiplicative semi-group $q^{B_2}$.  
Consider the Hecke algebra $\Hc=\langle (\Ga_2g\Ga_2) \rangle=\otimes_p\Hc_p$
 generated by all double coset classes $(\Ga_2g\Ga_2)$ with $g\in GSp_4(\Q)$.
Then we have that $\Hc_p=\Q[x_0^\pm, x_1^\pm, x_2^\pm]^{W_2}$ ($W_2$ the
 Weyl group)  and one has a $\Q$-algebras homomorphism $\la_f:\Hc\to \C$
 given by
 $f|X= \la_f(X)f$, $X\in \Hc$,  and $\al_j$ are defined as $\la_f(x_j)$,
 $j=0, 1, 2$.
In the notation of Andrianov, 
\begin{align}
Z_f(s)=L(s-k+(3/2), \pi_f, r_4) = \prod_p\det (1_4-h_p p^{k
-(3/2)})^{-1}p^{-s}
\end{align}
 is called the {\sl spinor $L$ function} of $f$, and he proved that it coincides
 with a linear combination of the Dirichlet series $\displaystyle{
 L(s, f, \xi_0)=\sum_{m=1}^\infty
\frac A_{m\xi_0}, m^s, }$
 where $\xi_0>0$ is a positive definite matrix of a fixed discriminant
 $-\det\xi_0$. 
Starting from this identity, he obtained an
 integral representation for $Z_f(s)$ using the group $GL_{2,K}$ where
$K=\Q(\sqrt{-\det\xi_0})$ an {\sl imaginary quadratic field}.
This integral representation implied an analytic continuation of $Z_f(s)$ to
 the whole complex plane and the functional equation of the type 
\begin{align}
\Psi_f(s)=\Ga_\C(s)\Ga_\C(s-k+2)Z_f(s)=(-1)^k\Psi_f(2k-2-s).
\end{align}
 where $\Ga_\C(s)=2(2\pi)^{-s}\Ga(s)$ is the standard $\Ga$-factor.
Its analytic properties were
 studied by A. N. Andrianov [And74] but still little is 
 known about  algebraic and arithmetic properties of critical values of this
 function; however, 
 the general Deligne conjecture on critical values of $L$-functions
 predicts that algebraicity properties could exist
 only for $s=k-1$ (see [Bo86, Fu-Sh, Ko-Ku] for evidences and discussions).

\medskip\noindent
The work of A.N.Andrianov was extended by I.I.Piatetski-Shapiro [PShBud],
 [PshPac]
 to arbitrary $F$ using a quadratic extension $K/F$ and the following
 construction.
Put 
$$
 V=K^2=\left\{x={x_1\choose x_2}, x_j\in K, j=1, 2\right \}
$$
 then $V$ may be viewed as a four dimensional $F$ vector space,
 $\dim_FV=4$, and define $\rho(x, y)=\Tr_{K/F}(x_1y_2-x_2y_1)$.
Let us consider the following $F$-algebraic group
\begin{align}
G=\{g\in GL_{2,K}\ |\ \det g\in GL_{1, F}\},\ \ SL_2(K)\subset G(F)\subset
GL_2(K)
\end{align}
 then there is an imbedding of $F$-algebraic groups
 $i:G\hookrightarrow G_\rho$
 because $x_1y_2-x_2y_1 = \det (x, y)$ and $\det (gx, gy)=\det g \cdot\det (x,
 y)$,  so that $\rho(gx, gy)=\det g\cdot\rho(x, y)$.
Note that $SL_2(\A_K)\subset G(\A_F)\subset GL_2(\A_K)$ and
 $G(\A_F)\hookrightarrow G_\rho(\A_F)=GSp_4(\A_F)$.
It turns out that there is an integral representation for $L(s, \pi, r)$ of the
 following type:
\begin{align}\label{(1.5)}
L(s, \pi, r) = 
\int_{G(F)C(\A_F)\bs G(\A_F)}\ph(i(g))E(g, s)dg := I_\pi (s) 
\end{align}
 where $\ph$ is an automorphic form on $G_\rho(\A_F)=GSp_4(\A_F)$ from the
 representation space of $\pi$, 
$C(\A_F)$ the center of $G(\A_F)\subset GL_2(\A_K)$, 
$E^\Phi (g, s)$ is a certain Eisenstein series on $G(\A_F)\subset
GL_2(\A_K)$ 
 attached to a Schwartz function $\Phi\in {\Sc}(V_\A)$
 ([PshPac], \S 5).

\Subsection{Initial idea of a $p$-adic construction.} 
Let $p\ge 5$ be a prime number.
We consider the
 case of two totally real fields $K\supset F$ and a representation $\pi_f$
 attached to a holomorphic Siegel-Hilbert cusp form $f(z)= \tilde \ph$
 of scalar weight $k=(k, \dots, k)$ on the Siegel-Hilbert half plane
\begin{align}
H_{2,F}=H_2\times \cdots \times H_2\ \ (n\ {\rm copies});
\end{align}
 in this case there is also a critical value $s=k-1$ for $L$-functions of the
 type $L(s, \pi_f, \otimes \chi, r)$ where $\chi$ is a  character of finite
 order of $\A_F^\times/F^\times$ . 
According to general conjectures on motivic $L$-functions there
 should exist  $p$-adic $L$-functions which interpolate $p$-adically their
 critical values, see [Co], [Co-PeRi], [PaIF].
However in our present construction instead of $p$-adic interpolation of 
 their special
 values of the type  $L(k-1, \pi_f \otimes \chi, r)$ we use directly a
 $p$-adic
 version of (\ref{(1.5)}) using techniques of $\La$-adic modular forms (see Section 3).
We hope that the resulting $p$-adic $L$-function provide also the above
 $p$-adic interpolation.

\Subsection{$\La$-adic modular forms.}
Let us consider the Iwasawa algebra [Iw] $\La=\Z_p[[T]]\cong \Z_p [[\Ga]]$
 as the completed group ring of the profinite group
 $\Ga=1+p\Z_p=\langle1+p\rangle \subset \Z_p^\times$.
We shall view elements of its quotient field $\Lc={\rm Quot}\La$ as
 $\C_p$-meromorphic functions with a finite number of poles
 on the unit disc $U = \{t\in \C_p| \ |\ |t|_p<1\}\subset \C_p$.
According to the theorem of Kubota-Leopoldt [Ku-Le],
 there exists a unique element
$g(T)\in \Lc$ such that for all $k\ge 1$, $k\equiv 0 \bmod (p-1)$
$$
g((1+p)^k-1)=\z^*(1-k)
$$
 where $\z^*(1-k)$ denotes the special value at $s=1-k$ of the Riemann
zeta-function with a modified Euler $p$-factor:
$\z^*(s)=(1-p^{-s})\z(s)$.
One could also start from positive values $s=k, k\equiv 0\bmod (p-1)$, 
and construct a $p$-adic zeta function $\z_{+,p}$ which interpolate 
$\displaystyle{k\mapsto
\z_+^*(k)=\frac \Ga(k), (2\pi i)^k, \z (k)(1-p^{k-1})}$ (see [Colm98]) and
 satisfies the following "functional equation" $\z_{+,p}(s)=2\z_p(1-s)$. 

\begin{defi}[The Serre ring] $\La [[q]]$ is the ring of all formal
 $q$-expansions with coefficients in $\La$:
$$
\La[[q]]=\{f=\sum_{n=0}^\infty a_n(T)q^n\ |\ a_n(T)\in \La\};
$$
\end{defi}

\begin{defi}
The $\La$-module $M(\La)\subset \La[[q]]$ of
{\sl $\La$-adic modular forms} (of some fixed level $N$, $(N, p)=1$ is
generated by  all $f=\sum_{n=0}^\infty a_n(T)q^n\in \La [[q]]$ such that for each
$k\ge 5$, $k \gg 0$
 the specialisation 
$$
f_k=f|_{T=(1+p)^k-1}\in \Z_p[[q]]
$$
is a classical modular form of weight $k$ and level $Np$.
In more precise terms $f$ is given by a $p$-adic measure $\mu_f$ on $\Z_p^\times$
with values in $\Z_p[[q]]$ such that the integrals 
\begin{align}
\int_{\Z_p^\times}x_p^k\mu_f=f_k
\end{align}
are classical modular forms.
\end{defi}

\begin{exam}
[The $\La$-adic Eisenstein series] $f\in M(\La)$ 
 (of level $N=1$) is defined by 
\begin{align}
f_k=\frac \z^*(1-k), 2,  + \sum_{n\ge 1}\si_{k-1}^*(n)q^n,
\ \ \ \si_{k-1}^*(n)=\sum_{d|n,p\nd d}d^{k-1}. 
\end{align}
\end{exam}

\begin{exam}
[Hida's families]
$f$ are elements of 
$$
S^{\rm ord}(\La)=
eS(\La),\ \ e=\lim_{n\to\infty}U_p^{n!}
$$
($U_p(\sum_{n\ge 0}a_nq^n)=\sum_{n\ge 0}a_{pn}q^n$ is the Atkin $U$-operator), 
$S(\La)$ is the $\La$-submodule of $\La$-adic cusp forms.
\end{exam}

\subsubsection*{The Hilbert modular case.} 
According to the classical theorem of
 Klingen [Kli], for a totally real field $K$ and for $k\ge 1$ the special values
 $\z_K(1-k)$ are rational numbers where $\z_K(s)$ is the Dedekind zeta function
 of $K$.

\noindent
{\sl The Deligne-Ribet $p$-adic zeta function} [De-Ri] interpolates $p$-adically
 these special values as an element $g_K\in \Lc$:
 for all positive
 integers
 $k\equiv 0 (\bmod (p-1))$, $g_K((1+p)^k-1)=\z^*_K(1-k)$, where
 $\z^*_K(s)=\prod _{{\mathfrak p}\ | p}(1-\Nc {\mathfrak p}^{-s})\z_K(s)$ is the Dedekind
zeta function of $K$ with all the ${\mathfrak p}$-factors over $p$
 removed from its Euler product.
The function  $\z_{K,p}(s)=g_K((1+p)^{1-s}-1)$ is
 analytic for all $s\in \Z_p\backslash 1$ with values in $\Q_p$ and it is called
the {\sl Deligne-Ribet $p$-adic zeta function}. 
It has the following properties: $\z_{K,p}(1-k)=\z^*_K(1-k)$ for all positive integers $k\equiv 0 (\bmod (p-1))$,
 and its residue $\Res_{s=1}\z_{K, p} (s)$ was computed by Colmez [Colm88]:
$\displaystyle{\Res_{s=1}\z_{K,p}=\frac 2^dh_KR_pE_p(1), w_K \sqrt{D_K},} $ where
 $d=[K:\Q]$, $E_p(s)=\prod _{{\mathfrak p}\ | p}(1-\Nc {\mathfrak p}^{-s})$, $R_p$ the
$p$-adic regulator of $K$ (which does not vanish according to the {\sl Leopoldt
conjecture}).

\noindent {\sl A $\La$-adic Hilbert modular form} could be defined
as a formal Fourier expansion  
$$
f= \sum_{0{\ll\atop =}\eta\in L_K}a_\eta q^\eta \in \La[[q^{L_F}]]
\ \ (L_K\subset K\  {\rm
a\ lattice})
$$
( $\eta$ runs over totally positive elements or 0) whose
 appropriate specialisations  are classical Hilbert modular form.
When  $h_K> 1$ one needs to consider collections of
 such series $\{f_\la\}\ (\la = 1, 2, \dots, h_K)$
 in order to be able to use the action of the Hecke algebra.
$\La$-adic Hilbert modular forms
 were used by Wiles in his proof of the Iwasawa
 conjecture over totally real fields (see [Wi90] where a precise definition of
 a $\La$-adic Hilbert modular form is contained in Section 3). 
It is required
 that
 for all appropriate sufficiently large $k$ the specialization
$f_k=f|_{T=(1+p)^k-1}$
 is the Fourier expansion of a classical Hilbert modular form.
As over $\Q$, the first natural example of a $\La$-adic Hilbert modular form is
 given by a $\La$-adic Eisenstein series (more precisely, this series is given
 by the Katz-Hilbert-Eisenstein measure, see [Ka78]). 
Also, Hida's theory could be
 extended to the Hilbert modular case and even to the general case
 of cohomological modular forms on $GL_{2, K}$ over an arbitrary number field
$K$ (see [Hi94]).

\subsection*{The Siegel-Hilbert modular case.} 
A $\La$-adic Siegel-Hilbert
 modular form could be defined as a formal Fourier expansion 
$$
f= \sum_{\xi\in L_{2,F}}A_\xi q^\xi \in \La[[q^{L_{2,F}}]]
\ \ (L_{2,F}\subset M_{2,F})\  
$$
 ($L_{2,F}$ is
 the semi-group of all symmetric totally  non-negative matrices $\xi$ in  a 
 sublattice of $M_{2,F}$) 
 whose
 appropriate specialisations  $f_k=f|_{T=(1+p)^k-1}$ are classical Siegel-Hilbert
 modular form. 
The first example of a $\La$-adic Siegel-Hilbert modular form is given
 by an Eisenstein series (for $F=\Q$ these series are described in [PaSE]). 
It seems that Hida's
 theory also could be extended to the Siegel-Hilbert modular case 
[Hi98], [Til-U],[Til].

\Subsection{$p$-adic $L$-functions.} Recall that we consider the
 case of two totally real fields $K\supset F$ and an irreducible representation
 $\pi=\pi_f$ attached to a holomorphic Siegel-Hilbert cusp form $f(z) = \tilde
\ph$ of scalar weight $k=(k, \dots, k)$ on the Siegel-Hilbert half plane
$$
H_{2,F}=H_2\times \cdots \times H_2\ \ (n\ {\rm copies});
$$
Then we rewrite the integral representation (1.5) in the form
 of the Petersson scalar product over $K= F(\sqrt D)$:
$$
I_\pi(1/2) = 
\langle \tilde i^*\tilde \ph,\tilde E(s, \mu)\rangle_K\eqno(4.1)
$$
 where $i$ denotes both the imbedding $i:G\hookrightarrow G_\rho$ and the
 corresponding modular imbedding
$$
i:H_F\times H_F \to H_{2,F},
\ \ H_{F}=H\times \cdots \times H\ ;
H_{2,F}=H_2\times \cdots \times H_2\ \ (n\ {\rm copies});\eqno(4.2)
$$
If we write this imbedding in coordinates it takes the form 
$$
(z_1, z_2)\mapsto Z(z_1, z_2) = C{\rm diag} \{z_1, z_2\}{}^tC\ (z_1, z_2\in H_F)
$$ 
 (see [Shi78], [Wi90], p. 521), where  we could take $C=\frac 1, 2, \mat 1, 1,
1/\sqrt
 D, -1/\sqrt D, \in M_2(F)$ so that $i^*\tilde\ph=\tilde \ph\circ i$ is a
  holomorphic Hilbert modular form with an explicitely given Fourier expansion.
 If $\tilde \ph = \sum_{\xi\in L_{2,F}}A_\xi q^\xi$ then
 $i^*\tilde\ph = \sum_{\eta\in L_F}a_\eta q^\eta$ where each Fourier coefficient 
 $a_\eta$
 is a finite sum of certain $A_\xi$: 
$$
a_\eta=\sum_{\xi: \eta=\frac 1,4, (\xi_{11}+\xi_{22}+2\xi_{12}/\sqrt D)}A_\xi,
$$ 
so that the map $\tilde \ph \mapsto 
 i^*\tilde \ph$ could be defined in terms of their formal $q$-expansion.
For the $\La$-adic construction let us take a $\La$-adic Siegel-Hilbert cusp 
 form
 $\tilde \ph$ on $GSp_{4, F}$ then $i^*\tilde \ph$ is a $\La$-adic Hilbert
 modular form over $K$ which is explicitely described as a formal
 Fourier expansion.
Now let us take $G$ to be the $\La$-adic Hilbert-Eisenstein series for
 $GL_{2, K}$.
In order to define the Petersson product 
$$
\langle \tilde i^*\ph, G \rangle_K\eqno(4.3)
$$
 we use the Eisenstein projection $1_{\rm Eis}(i^*\tilde\ph)$ (the projection in
 the $\Lc$-vector
 space $M(\Lc)$ to the (finite-dimensional) $\Lc$-subspace $Eis_K(\Lc)$ of
 Hilbert-Eisenstein series with an explicitely given base coming from the
Katz-Hilbert-Eisenstein $p$-adic measure). The projection $1_{\rm
Eis}(i^*\tilde\ph)$ could be explicitely computed
 using  the Fourier expansions of $i^*\tilde\ph$ and 
 of the Fourier expansions of a
 $\Lc$-basis of $Eis_K(\Lc)$.
$$
\langle \tilde i^*\ph, G \rangle_K=\langle 1_{\rm Eis}(i^*\tilde\ph),
G\rangle_K. 
$$
Then  we are reduced to the case of $\langle G_1, G_2\rangle_K$, where $G_1$ 
 and $G_2$ are two normalized Hilbert-Eisenstein series, and in order to
 define their Petersson product we use the method of Rankin-Selberg.

\medskip\noindent 
Let us recall a classical formula 
$$
(f,g)=\frac \pi, 3, \frac \Ga(k), (4\pi)^{k}, \Res_{s=k}L_{f,g}(s)
$$
for the Petersson product
$\displaystyle{
(f,g)=\int_{\Ga\backslash \H}f(z)\overline{g(z)}y^{k-2}\, dx\,dy}$
 (see [Ra39, Za81]) where  

\noindent $L_{f,g}(s) = \sum_{n=1}^\infty a_n\overline{b}_n
n^{-s}$ denotes the Rankin $L$-function of two
 holomorphic modular forms of weight $k$ on $SL_2(\Z)$ ,
with at least one of them a cusp form (i.e. $a_0b_0=0$ ):
$f(z) = \sum_{n=0}^\infty a_ne^{2\pi i nz}$ and $g(z) = \sum_{n=0}^\infty
b_ne^{2\pi i nz}$.  
This equality makes it possible to define the
Petersson scalar product (a {\sl renormalized value}) $(G_k, G_k)$ where
$\displaystyle{
G_k=-\frac B_k, 2k, + \sum_{n=1}^\infty \si_{k-1}(n)}$ ($k\ge 4, k$
even). We have [Za81, p.435]: 
$$
L_{G_k,G_k}(s) = \sum_{n=1}^\infty
\si_{k-1}(n)\si_{k-1}(n) n^{-s} = \frac
{\z(s)\z(s-k+1)^2\z(s-2k+2)},{\z(2s-2k+2)}, 
$$
which implies
\begin{align*}
(G_k, G_k)&=(-1)^{k/2-1}\frac \Ga(k)\Ga(k-1), 2^{3k-3}\pi^{2k-1},\z(k)\z(k-1)\\ &=
i^{3k-3}2^{2-k}\frac \Ga(k), (2\pi i)^{k}, \z (k) \frac \Ga(k-1), (2\pi
i)^{k-1}, \z (k-1) . 
\end{align*}

\medskip\noindent 
We see that if $G_1, G_2$ were two cusp forms of weight $k$ their Petersson
product would
 essentially coincide with a normalized residue of the Rankin zeta function
 $L_{G_1, G_2}(s)$ at $s=k$.
In the case of normalised Eisenstein series the
 Rankin zeta function $L_{G_1, G_2}(s)$ is explicitely evaluated via Rankin's
 lemma as a product of abelian Dirichlet $L$-functions.
Let now $G_1 = \{G_{1,k}\}$, $G_{2,k}= \{G_{2,k}\}$ denote two $p$-adic
families of Hilbert Eisenstein series. 
We may define the $I_{G_1, G_2}=\langle G_1, G_2\rangle_K$ as
an element of $\Lc$ such that for all $k \gg 0$ 
$$
(G_1, G_2)=\Res_{s=k}L_{G_{1,k}, G_{2,k}}(s)\in \Q_p \ (s\in \Z_p)
$$
 in a similar way as in [Za81] and [Ko-Za] as the normalised $p$-adic residue
 of the $p$-adic Rankin convolution
 $L_{G_1, G_2}(s)$ (which is defined in terms of the corresponding Deligne-Ribet
 $p$-adic zeta function).

\begin{main}
 Let $\tilde\ph$ be a $\La$-adic Siegel-Hilbert modular eigenform
 then there exists a canonically defined  element 
$$
I_{\tilde \ph, p}=\langle \tilde i^*\ph, G \rangle_K\in \Lc
$$
$i^*\tilde \ph$ the $\Lambda$-adic pullback of $\ph$,   
$i^*\tilde \ph$ is a $\La$-adic Hilbert
 modular form over $K$ explicitely described by its Fourier expansion, 
 $G$ is a certain  $\La$-adic Hilbert-Eisenstein series,
such that the function $I_{\tilde \ph, p}$
 gives a $p$-adic interpolation of
 the  residue of the normalized $p$-adic Rankin $L$ function 
$ L^*_{i^*\tilde\ph_k, G_k}(s)$ (at $s=k$),
 the scalar weight of a specialisation $\tilde\ph_k)$:
$$
I_{\tilde \ph, p}|_{T=(p+1)^k-1}=\Res_{s=k}L^*_{i^*\tilde\ph_k, G_k}(s)\ (s\in
\Z_p) 
$$
\end{main} 

\Subsection{$p$-adic families of automorphic representations.}
We use the occasion to discuss here the following general definition
of a {\sl $p$-adic family of automorphic representations} (or of a {\sl 
$\La$-adic automorphic form}). 
We shall view the Iwasawa algebra $\La$ as the algebra ${\rm Meas}(\Z_p, \Z_p)$
of all $\Z_p$-measures on $\Z_p$ (with the additive convolution
 as a  multiplication).
Let $V_\Q\subset C(G(\A_F))$ be a certain $\Q$-vector space
 of (complex-valued) continuous functions on the adelic group $G(\A_F)$ over
 a number field $F$.
We suppose that $V_\Q$ has an integral structure $V_\Z\subset V_\Q$ so that
 $V_\Q=V_\Z\otimes \Q$. 
Put $V_p=V_\Z\hat \otimes \Z_p$ (the completed tensor product).
Define $D_p(V_p)={\rm Meas}(\Z_p, V_p)$ (as a module over $\La=D_p(\Z_p)$).
 
\begin{defi}
A $p$-adic family of automorphic representations
 on $G$  is a $p$-adic measure $\ph \in D_p(V_p)$
 such that for almost all positive integers $k$ we have that
 the integral $\int_{\Z_p}x^k\ph=\ph_k\in V_p$ belongs to $V_\Z$
 and the function $\ph_k$ generates an automorphic representation 
 $\pi_k$ of $G(\A_F)$. 
We call such $\ph$ a $\La$-adic automorphic form on $G(\A_F)$.

Let $AF_G(\La)$ denote the $\La$-module generated by such elements $\ph$.
An element $\ph$ is called an eigenform if the representations $\pi_k$ are all
 {\sf irreducible}.
\end{defi}
 
 \medskip\noindent 
A natural example of such a vector space $V$ for the group   
$\GL_2$ over $\Q$ comes from 
  holomorphic functions $f=\sum_{n=0}^\infty a_n \exp (2\pi i n z)$
  having rational Fourier coefficients $a_n\in \Q$ with bounded denominators,
  i.e. for which there exists a positive integer $N=N(f)$ such that $Na_n\in
\Z$. 
However there are other ways to attach such a vector space $V$ to $G$
 by considering cohomology groups of the corresponding locally-symmetric 
 spaces and automorphic forms $\ph$ on $G(\A_F)$ represented by rational
 cohomology classes ([Ko-Za]). 
Put $AF_G(\Lc)=AF_G(\La)\otimes \Lc$.
We hope that one could find in this way a general
 construction of
 $p$-adic automorphic $L$ functions $L_{\pi, r, p}$ as certain $\Lc$-linear forms
 $l=l_{G, r}$ on the $\Lc$-vector space
 $AF_G(\Lc)$.
Such a linear form should play a role of an integral representation for 
 $p$-adic $L$-functions: $L_{\pi_k, r, p}=l_{G, r}(\ph)|_{T=(1+p)^k-1}$.
A natural example of such a linear form comes from the $\La$-adic Petersson
 product of Hida which provides a construction of $p$-adic $L$-functions for
 $GL_2\times GL_2$ [Hi91].

On the other hand, there exist nice constructions of $p$-adic families of
 Galois representations attached to automorphic forms
 ($\La$-adic Galois representations, see [Hi86], [Til-U]) which played
 an important role in the work of Wiles [Wi95].
It would be interesting to formulate a general $\La$-adic Langlands conjecture
 relating $\La$-adic automorphic forms and $\La$-adic Galois representations.

\bigskip
\noindent
{\large\bf  References for the Appendix}
\medskip

\par
\medskip
\newdimen\brawidth\setbox0=\hbox{[PaAdm1]\ }\brawidth=\wd0

\def\ref#1,#2!{\smallskip\noindent
\hangindent=\brawidth
\hangafter=1
\hbox to\brawidth{[#1]\hfill}{\it #2}:\ }

\ref AndBud,Andrianov A.N.!Zeta-functions and the Siegel modular
 forms. Proc. Summer School Bolya-Janos Math.Soc (Budapest, 1970), Halsted,
 N.-Y., 1975, 9-20

\ref And74,Andrianov A.N.!Euler   products    attached \ to \ Siegel
 modular  forms \ of \ degree 2 , Uspekhi Mat. Nauk 29 (1974) 44-109
(in Russian)

\ref Bo86,B\"ocherer S.! Bemerkungen \"uber
 die Dirichletreihen von Koecher and Maass. Math. Gottingensis, 
 Schriftenr. d. Sonderforschungsbereichs Geom.Anal. 68 (1986)

\ref Bo-Sch,B\"ocherer S.,Schmidt C.-G.!$p$-adic measures attached to Siegel
modular forms (to appear in Annales de l'Institut Fourier)

\ref Co,Coates J.!On $p$-adic $L$-functions. S\'eminaire Bourbaki, 40eme annee,
1987-88, n${}^\circ$ 701, Asterisque (1989) 177-178.

\ref Co-PeRi,Coates J., Perrin--Riou B.!On $p$-adic $L$-functions
attached to motives over ${\bf Q}$. Advanced Studies in Pure Math. 17
(1989), 23--54

\ref Colm88,Colmez, P.! 
R\'esidu en $s=1$ des fonctions z\^eta $p$-adiques. 
 Invent. Math. 91 (1988), no. 2, 371-389

\ref Colm98,Colmez, P.!Fonctions $L$ $p$-adiques. S\'eminaire Bourbaki, 
51 \`eme ann\'ee, 1998-99, n${}^\circ$ 851. Novembre 1998

\ref  De79,Deligne P.!Valeurs de fonctions $L$ et p\'eriodes d'int\'egrales.
\ Proc. Symp. Pure Math AMS 33 (part 2) (1979) 313 - 342

\ref  De-Ri,Deligne P.,Ribet ,K.A.!Values of Abelian $L$-functions at
negative integers over totally real fields. Invent. Math. 59
(1980) 227-286

\ref Fu-Sh, Furusawa, Masaaki; Shalika, Joseph A.!The fundamental lemma
 for the Bessel and Novodvorsky subgroups of ${\rm
GSp}(4)$. C. R. Acad. Sci. Paris S\'er. I Math. 328 (1999), no. 2, 105--110

\ref Ge-PSh,Gelbart S.,Piatetski-Shapiro I.I., Rallis S.!Explicit
constructions of automorphic $L$ - functions. Springer-Verlag, Lect. Notes in
Math. N 1254 (1987) 152p.

\ref  Hi86,Hida H.!Galois representations into $GL_{2}({\bf Z}_{p}[[X]])$ 
 attached to ordinary cusp forms, Invent. Math. 85 (1986) 545-613
 
\ref Hi90,Hida, H.!Le produit de Petersson et de Rankin $p$-adique. 
 S\'eminaire de Th\'eorie des Nombres, Paris 1988--1989, 87-102,
 Progr. Math., 91, Birkh\"auser Boston, Boston, MA, 1990 

\ref  Hi91,Hida H.!On $p$-adic $L$-functions of $GL(2)\times GL(2)$ over
 totally real fields, Ann. l'Inst. Fourier, 40, no.2 (1991) 311-391

\ref Hi94,Hida, H.! 
$p$-adic ordinary Hecke algebras for ${\rm GL}(2)$. 
Ann. Inst. Fourier (Grenoble) 44 (1994), no. 5, 1289-1322

\ref Iw,Iwasawa K.!Lectures on $p$-adic $L$-functions, Ann. of Math.
 Studies, N 74. Princeton University Press, 1972

\ref Ka78,Katz, N. M.! 
$p$-adic $L$-functions for CM fields. 
Invent. Math. 49 (1978), no. 3, 199-297

\ref Kli,Klingen, H.!\"Uber die Werte
 der Dedekindschen Zetafunktion. Math. Ann. 145 1961/1962 265-272

\ref Ku-Le,Kubota T., Leopoldt H.-W!Eine $p$-adische Theorie der
Zetawerte. J. reine angew. math. 214/215 (1964) 328-339

\ref Ko-Ku,Kohnen, W.,  Kuss, M! Some numerical computations concerning
 spinor zeta functions in genus 2 at the central point. Preprint, 1999

\ref Ko-Za,Kohnen W., Zagier D.!Modular forms with rational periods. In: Modul.
Forms Symp., Durham, 30 June - 10 July 1983. Chichester, 1984 197-249
 
\ref Lau,Laumon G.!Sur la cohomologie \`a supports compacts des vari\'et\'es
de Shimura pour ${\rm GSP}(4)_\Q$. Comp. Math.,105 (1996) 267-359.

\ref Man,Manin Yu.I.!Non-Archimedean integration and
 Jacquet-Langlands $p$-adic $L$-functions.
 Russ. Math. Surveys 31, N 1 (1976) 5-57

\ref PaLNM,Panchishkin A.A.!Non-Archimedean $L$-functions of
 Siegel and Hilbert modular forms, Lecture Notes in Math., 1471, 
 Springer-Verlag, 1991, 166p.

\ref PaAdm,Panchishkin A.A.!Admissible Non-Archimedean standard zeta functions
 of Siegel modular forms, 
 Proceedings of the Joint AMS Summer Conference on Motives,
 Seattle, July 20--August 2 1991, Seattle, Providence, R.I., 1994, vol.2, 
 251 -- 292 

\ref PaIF,Panchishkin A.A.!Motives over totally real fields and $p$-adic
 $L$-functions,  Annales de
 l'Institut Fourier, Grenoble, 44 (1994) 989-1023

\ref PaSE,Panchishkin A.A.!On the Siegel--Eisenstein measure and
 its applications.
(to appear in the Israel Journal of Mathemetics, 1999)

\medskip
\noindent
\ref  PaViet,Panchishkin A.A.!Non-Archimedean Mellin transform and
 $p$-adic $L$ Functions.
Vietnam Journal of Mathematics, 1997, N3, 179-202 

\ref PShBud,Piatetski-Shapiro I.I.!Euler subgroups, Proc.
 Summer School Bolya-Janos Math.Soc (Budapest, 1970), Halsted, N.-Y.,
 1975, 597-620 

\ref PShPac,Piatetski-Shapiro I.I.!$L$-functions for $GSp_4$. Pacific J.Math,
Olga Tausski-Todd memorial issue, (1998) 259-275

\ref PSh-R,Piatetski-Shapiro I.I., Rallis S.!$L$-functions of
automorphic forms on simple classical groups. In: Modul. Forms Symp.,
Durham, 30 June - 10 July 1983. Chichester, 1984, 251-261

\ref Ran39,Rankin R.A.! Contribution to the theory of Ramanujan's
function $\tau (n)$ and similar arithmetical functions. I.II.Proc.
Camb. Phil. Soc 35 (1939) 351-372

\ref Ran52,Rankin R.A.!The scalar product of modular forms. Proc.
 London math. Soc. 2 (1952) 198-217

\ref Shi76,Shimura G.!The special values of the zeta functions associated
with cusp forms. Comm. Pure Appl. Math. 29 (1976) 783-804

\ref Shi77,Shimura G.!On the periods of modular forms. Math. Annalen 229
(1977) 211-221

\ref Shi78,Shimura G.!The special values of zeta functions associated
with Hilbert modular forms. Duke Math. J. 45 (1978) 637-679

\ref Shi83,Shimura G.!On Eisenstein series. Duke Math. J. 50 (1983)
417-476

\ref Tay,Taylor R.!On the $l$-adic cohomology of Siegel threefold. Invent.
 Math. 114 (1993) 289-310

\ref Til-U,Tilouine J., Urban E.!Several variable $p$--adic
 families of Siegel-Hilbert cusp eigenforms and their Galois representations.
 Ann. scient. \'Ec. Norm. Sup. 4${}^{\rm e}$ s\'erie, 32 (1999) 499-574

\ref Til,Tilouine J.!Deformations  of Siegel-Hilbert Hecke
 eigensystems and their Galois representations. Contemporary Math. 210 (1998)
195-225

\ref Wi88,Wiles A.!On ordinary $\lambda$-adic representations
 associated to modular forms.
 Invent. Math. 94, No.3, 529-573 (1988)

\ref Wi90,Wiles A.!The Iwasawa conjecture for totally real fields. 
Ann. of Math. (2) 131 (1990), no. 3, 493--540. 

\ref Wi95,Wiles A.!
Modular elliptic curves and Fermat's Last Theorem.
Ann. Math., II. Ser. 141, No.3 (1995) 443-55

\ref Za81,Zagier, Don! The Rankin - Selberg method for automorphic
 functions which are not of rapid decay. J. Fac. Sci. Univ.
 Tokyo Sect. IA Math. 28 (1981), no. 3, 415-437 (1982).


\bibliographystyle{plain}

\begin{thebibliography}{}

\end{thebibliography}


\begin{thebibliography}{10xxxx}

\bibitem[Am-V]{Am-V}
{\sc Amice,  Y.} and {\sc Vélu,  J.}, 
{\em Distributions $p$-adiques associ\'ees aux
s\'eries de Hecke}, 
 Journ\'ees Arithm\'etiques de Bordeaux (Conf. Univ.
Bordeaux, 1974), Ast\'e\-risque no. 24/25, Soc. Math. France, Paris
1975, 119 - 131


\bibitem[Ar89]{Ar89}
{\sc Arthur}, J., 
{\em Unipotent automorphic representations: Conjectures.} in:
 Orbites unipotentes et représentations, II, Astérisque 171-172, Soc. Math. France,
Montrouge, 1989, 13-71. 

\bibitem[Boe84]{Boe84}
{\sc B\"ocherer}, S., \ {\em \"Uber die Fourierkoeffizienten Siegelscher Eisensteinreihen}, Manuscripta
Math., 45 (1984), 273-288.



\bibitem[Boe85]{Boe85}
{\sc B\"ocherer}, S., \ {\em \"Uber die Funk\-tio\-nal\-glei\-chung 
auto\-mor\-pher $L$--Funk\-tio\-nen zur Sie\-gel\-scher Modul\-gruppe.} 
J. reine angew. Math. 362 (1985) 146--168

\bibitem[Boe-Pa9]{Boe-Pa9}
{\sc B\"ocherer}, S.,  {\sc Panchishkin},  A.A.,
{\it $p$-adic Interpolation of Triple $L$-functions: Analytic Aspects}. In:
Automorphic Forms and $L$-functions II: Local Aspects -- David Ginzburg, 
Erez Lapid, 
and David Soudry, 
Editors, AMS, BIU,  2009, 313 pp.; pp.1-41

\bibitem[Boe-Pa11]{Boe-Pa11}
{\sc B\"ocherer}, S.,  {\sc Panchishkin},  A.A.,
{\it Higher Twists and Higher Gauss Sums}
Vietnam Journal of Mathematics 39:3 (2011) 309-326


\bibitem[Boe-Schm]{Boe-Sch}
{\sc B\"ocherer, S.},  and {\sc Schmidt, C.-G.}, 
{\em $p$-adic measures attached to Siegel modular forms}, 
Ann. Inst. Fourier 50, \Numero 5, 1375-1443 (2000).

\bibitem[Co]{Co}
{\sc Coates, J.} {\em
On $p$--adic $L$--functions.} Sem. Bourbaki, 40eme annee,
1987-88, n${}^\circ$ 701, Asterisque (1989) 177--178.

\bibitem[Co-PeRi]{Co-PeRi}
{\sc Coates, J.} and {\sc Perrin-Riou, B.}, 
{\em On $p$-adic $L$-functions attached to 
motives over ${\Q}$},
Advanced Studies in Pure Math. 17, 23--54 (1989)

\bibitem[CGS]{CGS}
{\sc Jim Cogdell, Steve Gelbart, and Peter Sarnak, Coordinating Editors}, 
{\em
Ilya Piatetski-Shapiro,
In Memoriam.}
November 2010, Notices of the AMS, Volume 57, Number 10  p.1260-1275

\bibitem[ConSl98]{ConSl98}
{\sc Conway, J.H.}, {\sc Sloane, N.J.A.},
  Sphere Packings, Lattices and
Groups, Springer-Verlag, NY, 3rd edition, 1998.



\bibitem[CourPa]{CourPa}
{\sc Courtieu,}M.,  {\sc Panchishkin },A.A.,
{\em Non-Archimedean $L$-Functions and Arithmetical Siegel Modular Forms},
Lecture Notes in Mathematics 1471, Springer-Verlag, 2004 (2nd augmented ed.)

\bibitem[GPS]{GPS}
 {\sc Gelbart,} S.,  {\sc Panchishkin,} A., and  {\sc Shahidi,} S., {\it The $p$-Adic Eisenstein Measure and Fourier coefficients for SL(2)}
   Arxiv NT, 1001.1913 (2010) 


\bibitem[GeSha]{GeSha}
 {\sc Gelbart,} S.,  and  {\sc Shahidi}, F.
{\em Analytic Properties of Automorphic $L$-functions}, Academic Press, New York, 1988.

\bibitem[GRPS]{GRPS}  {\sc Gelbart S.,Piatetski-Shapiro I.I., Rallis S.}
{\it Explicit
constructions of automorphic $L$ - functions.}
 Springer-Verlag, Lect. Notes in Math.
N 1254 (1987) 152p.

\bibitem[Gue]{Gue}
 {\sc Guerzhoy}, P.
{\em On $p$-adic families of Siegel cusp forms in the Maass Spezialschaar.}
Journal für die reine und angewandte Mathematik 523 (2000),  103-112 
 
\bibitem[Ha81]{Ha81}
{\sc Harris, } M., {\em
The rationality of holomorphic Eisenstein series}, Inv. Math. 63 (1981), 305-310. 

\bibitem[Ha84]{Ha84}
{\sc Harris, } M., {\em Eisenstein Series on Shimura Varieties.}
Ann. Math., 119 (1984), No. 1  59-94



\bibitem[HLiSk]{HLiSk}
{\sc Harris, }M., {\sc Li,} Jian-Shu., {\sc Skinner,} Ch.M.,
{\em $p$-adic $L$-functions for unitary Shimura varieties.}
Documenta Math.
Extra volume~: John H.Coates' Sixtieth Birthday (2006), p.393-464

\bibitem[He27]{He27}
 {\sc Hecke, E.},
        {\it Theorie der Eisensteinschen Reihen und ihre Anwebdung auf
             Funktionnentheorie und Arithmetik},
        Abh. Math. Sem. Hamburg 5 (1927), p. 199-224.

\bibitem[Hi93]{Hi4}
 {\sc Hida, H.}, {\it Elementary theory of $L$-functions and Eisenstein series.}
 London Mathematical Society Student Texts. 26 Cambridge, 1993



\bibitem[Ike01]{Ike01}
{\sc Ikeda, T.},
{\em On the lifting of elliptic cusp forms to Siegel cusp forms of degree $2n$}, Ann. of
Math. (2) 154 (2001), 641-681.


\bibitem[Ike06]{Ike06}
{\sc Ikeda, T.},
{\em Pullback of the lifting of elliptic cusp forms and Miyawaki's Conjecture
}
Duke Mathematical Journal, \textbf{131}, 469-497 (2006)




\bibitem[Ka76]{Ka76}  
{\sc  Katz, N.M.}, {\em $p$-adic interpolation of real analytic Eisenstein
series.} Ann. of Math. 104 (1976) 459--571




\bibitem[Kawa]{Kawa}  
{\sc  Kawamura},Hisa-Aki, {\em On certain constructions of $p$-adic families of Siegel modular forms of even genus} 
ArXiv, 1011.6042v1 



\bibitem[KMS2000]{KMS2000} D. {\sc Kazhdan}, B. {\sc Mazur},
C.-G. {\sc Schmidt},
        {\it Relative modular symbols and Rankin-Selberg convolutions}. 
J. Reine Angew. Math. 519, 97-141 (2000). 

\bibitem[Kl62]{Kl62}  
{\sc Klingen H.}, {\em \" Uber die Werte der Dedekindschen Zetafunktionen}.
Math. Ann. 145 (1962) 265--272




\bibitem[Kl67]{Kl67}  
{\sc Klingen H.}, {\em Zum Darstellungssatz für Siegelsche Modulformen}.
 Math. Z. 102 (1967) 30--43


\bibitem[Kub]{Kub}  
{\sc Kubota}, T.,  {\em Elementary Theory of Eisenstein
Series}, Kodansha Ltd. and John Wiley and Sons (Halsted Press),
1973

\bibitem[LangMF]{LangMF}
{\sc
Lang}, Serge. {\em Introduction to modular forms. With appendixes by D. Zagier and Walter Feit.} 
Springer-Verlag, Berlin, 1995

\bibitem[Maa71]{Maa71} 
{\sc Maass H.}, {\em Siegel's modular forms and Dirichlet series}
Springer-Verlag, Lect. Notes in Math. N 216 (1971)


\bibitem[Miy89]{Miy89} 
{\sc Miyake,~Toshitsune},  
        {\em Modular forms}. Transl. from the Japanese by Yoshitaka Maeda.,
       Berlin etc.: Springer-Verlag. viii, 335 p. (1989).


\bibitem[Mi92]{Mi92}
{\sc Miyawaki, Isao},
{\em Numerical examples of Siegel cusp forms of degree 3 and their zeta-functions},
Memoirs of the Faculty of Science, Kyushu University, Ser. A, Vol. 46, No. 2 (1992), pp. 307--339.

\bibitem[Pa81]{Pa81}
{\sc Panchishkin,  A.A.}, 
{\em Complex valued measures attached to Euler products},  Trudy Sem. 
Petrovskogo 7 (1981) 239-244 (in Russian)

\bibitem[Pa82]{Pa82}
{\sc Pan\v ci\v skin} A.A.,  
 Le prolongement $p$-adique analytique de fonctions $L$ de 
 Rankin I,II.  C. R. Acad. Sci. Paris  294 (1982) 51-53, 227-230.

\bibitem[Pa91]{Pa91}
{\sc Panchishkin,  A.A.}, 
{\em Non--Archimedean $L$-functions of Siegel
and Hilbert modular forms}, Lecture Notes in Math., {\bf 1471}, 
Springer--Verlag, 1991, 166p.

\bibitem[Pa94]{Pa94}
{\sc Panchishkin,  A.A.}, 
{\em Admissible Non-Archimedean standard zeta functions
 of Siegel modular forms}, 
 Proceedings of the Joint AMS Summer Conference on Motives,
 Seattle, July 20--August 2 1991, Seattle, Providence, R.I., 1994, vol.2, 
 251 -- 292 
 
 \bibitem[PaSE]{PaSE}  {\sc Panchishkin, A.A.},
        {\it On the Siegel-Eisenstein measure and
its applications},
        Israel Journal of
Mathemetics, 120, Part B (2000) 467-509.

\bibitem[PaMMJ]{PaMMJ} {\sc Panchishkin,   A.A.},
{\em A new method of constructing $p$-adic $L$-functions associated 
with modular forms}, 
Moscow Mathematical Journal, 2 (2002), Number 
2, 1-16


\bibitem[PaTV]{PaTV} {\sc Panchishkin,   A.A.},
{\em Two variable $p$-adic $L$ functions attached to eigenfamilies of positive slope},
 Invent. Math. v. 154, N3 (2003), pp. 551 - 615


\bibitem[PaIAS]{PaIAS}
{\sc Panchishkin,  A.A.}, 
{\em On $p$-adic integration in spaces of modular forms
 and its applications},  
J. Math. Sci., New York 115, No.3, 2357-2377 (2003). 

 
\bibitem[PaMZ10]{PaMZ10}
{\sc Panchishkin}, A.A.,  
{\em Two modularity lifting conjectures for families of Siegel modular forms},  
Mathematical Notes 
Volume 88 (2010) Numbers 3-4, 544-551

 \bibitem[PaIsr11]{PaIsr11}  {\sc Panchishkin, A.A.},
        {\it 
Families of Siegel modular forms, L-functions
and modularity lifting conjectures.}
Israel Journal of
Mathemetics, 185 (2011), 343-368

\bibitem[PS1]{PS1}
{\sc Piatetski-Shapiro, I. I.}, {\em On the Saito-Kurokawa lifting}, Invent. Math. 71 (1983), 
309-338.


\bibitem[PS2]{PS2}
{\sc Piatetski-Shapiro,  I. I.}, {\em Some example of automorphic forms on Sp4}, Duke Math. J. 50 (1983), 55-106.

\bibitem[PS3]{PS3}
{\sc Piatetski-Shapiro, I. I.}, {\em L-functions for GSp4}, Olga Taussky-Todd: in memoriam, Pac.
J. Math. 181(3) (1997), 259-275.

\bibitem[PSR]{PSR}
{\sc Piatetski-Shapiro, I. I.,}  and {\sc Rallis, S.}, {\em A new way to get Euler products}, J. Reine
Angew. Math. 392 (1988), 110-124.


\bibitem[Si35]{Si35}
{\sc Siegel}, C. L., 
{\em   Über die analytische Theorie der quadratischen
Formen}. Ann. of Math. 36, 527-606 (1935)


\bibitem[Si39]{Si39}
{\sc Siegel}, C. L., 
{\em   Einf\" uhrung in die Theorie der Modulfunktionen
$n$-ten Grades.} Math. Ann. 116 (1939) 617-657

\bibitem[Si64a]{Si64a}
{\sc Siegel}, C. L.,  {\em Zu zwei Bemerkungen Kummers.} Nachrichten der Akademie der Wissenschaften in Göttingen. Mathematisch-physikalische Klasse, 1964, Nr. 6, 51-57 


\bibitem[Si64b]{Si64b}
{\sc Siegel}, C. L., 
{\em   Über die Fourierschen Koeffizienten
der Eisensteinsehen Reihen}. Mat. Fys. Medd. Danske
Mid. Selsk. 34, Nr.6 (1964)



\bibitem[Si70]{Si70}
{\sc Siegel}, C. L., 
{\em   \" Uber  die  Fourierschen Koeffizienten
von Modulformen.} Nachr. Acad. Wiss. G\" ottingen. II. Math.- Phys. Kl. 3
(1970) 15--56




 \bibitem[Sha81]{Sha81}
{\sc Shahidi, F.}, 
{\em On Certain L-functions}, Amer. J. Math, 1981.


 \bibitem[Sha88]{Sha88}
{\sc Shahidi, F.}, 
{\em On the Ramanujan Conjecture and Finiteness of Poles for
Certain L-functions}, Annals of Math., 1988.


\bibitem[Shi95]{Shi95}
{\sc Shimura} G., 
{\em Eisenstein series and zeta functions on symplectic
 groups}, Inventiones Math. 119 (1995) 539--584




\bibitem[MC]{MC}
{\sc Skinner},  C.  and {\sc Urban}, E.  {\em
The Iwasawa Main Cconjecture for  GL(2)}.
{\tt http://www.math.jussieu.fr/~urban/eurp/MC.pdf}

\bibitem[St81]{St81}
{\sc Sturm} J., 
{\it The critical values of zeta-functions associated to
the symplectic group.} Duke Math. J. 48 (1981) 327-350

\bibitem[Wa82]{Wa82}
{\sc Washington,  L.},
{\em Introduction to cyclotomic fields},
 Springer
 Verlag: N.Y. e.a., 1982



\end{thebibliography}

\end{document}